\magnification=1100
\input amssym
\input xy
\xyoption{all}
\overfullrule=0pt
\def\Der{{\rm Der}}
\def\Q{{\Bbb Q}}
\def\t{{\frak{t}}}
\def\r{{\frak{r}}}
\def\grt{{\frak{grt}}}
\def\u{{\frak{u}}}
\def\pls{{\frak{pls}}}
\def\ds{{\frak{ds}}}
\def\nz{{\frak{nz}}}
\def\dmr{{\frak{dmr}}}
\def\ls{{\frak{ls}}}
\font\wee=cmr8
\font\weeit=cmti8
\font\weebf=cmbx8
\centerline{\bf Elliptic double shuffle, Grothendieck-Teichm\"uller and mould theory}
\vskip .3cm
\centerline{Leila Schneps}
\vskip .9cm
{\narrower{\wee{\weebf Abstract.} In this article\footnote{}{\it 2010 Mathematics Subject Classification: 11M32}
we define an {\weeit elliptic double shuffle Lie
algebra} $\scriptstyle{\ds_{ell}}$ that generalizes the well-known {\weeit double shuffle
Lie algebra} $\scriptstyle{\ds}$ to the elliptic situation. The double shuffle, or
dimorphic, relations satisfied by elements of the Lie algebra $\scriptstyle{\ds}$ 
express two families of algebraic relations between multiple zeta 
values that conjecturally generate all relations.
In analogy with this, elements of the elliptic double shuffle Lie algebra 
$\scriptstyle{\ds_{ell}}$ are Lie polynomials having a dimorphic property called 
$\scriptstyle{\Delta}$-bialternality that conjecturally describes the (dual of the) set of 
algebraic relations between {\weeit elliptic multiple zeta values}, which arise
as coefficients of a certain elliptic generating series (constructed
explicitly in [LMS]) and closely related to the elliptic associator defined by 
Enriquez ([En1]).  We show that one of 
Ecalle's major results in mould theory can be reinterpreted
as yielding the existence of an injective Lie algebra morphism 
$\scriptstyle{\ds\rightarrow \ds_{ell}}$.  Our main result is the compatibility of this 
map with the tangential-base-point section $\scriptstyle{{\rm Lie}\,\pi_1(MTM)\rightarrow 
{\rm Lie}\,\pi_1(MEM)}$ constructed by Hain and Matsumoto and with
the section $\scriptstyle{\grt\rightarrow \grt_{ell}}$ mapping
the Grothendieck-Teichm\"uller Lie algebra $\scriptstyle{\grt}$ into the elliptic 
Grothendieck-Teichm\"uller Lie algebra $\scriptstyle{\grt_{ell}}$ constructed by Enriquez.
This compatibility is expressed by the commutativity of the following diagram
(excluding the dotted arrow, which is conjectural).}\par}
$$\xymatrix{
{\rm Lie}\,\pi_1(MTM)\ar@{^{(}->}[r]^{\ \ \ \ \rm Brown}\ar[d]_{\rm Hain}^{\rm Matsumoto}&\grt\ar[d]^{\rm Enriquez}
\ar@{^{(}->}[r]^{{\rm Furusho}}&\ds\ar[d]^{\rm Ecalle}\\
{\rm Lie}\,\pi_1(MEM)\ar[r]\ar[dr]&\grt_{ell}\ar@{.>}[r]\ar[d]&\ds_{ell}\ar[dl]\\
&\Der\,{\rm Lie}[a,b]& }\eqno(A)$$
\vskip .3cm
{\narrower{{\weebf R\'esum\'e.} {\wee Dans cet article, nous d\'efinissons une
{\weeit alg\`ebre de Lie de double m\'elange elliptique}
$\scriptstyle{\ds_{ell}}$ qui g\'en\'eralise l'alg\`ebre de Lie 
bien connue {\weeit de double m\'elange} $\scriptstyle{\ds}$ au cas 
elliptique. Les relations de double m\'elange (ou dimorphiques) satisfaites
par les \'el\'ements de l'alg\`ebre de Lie $\scriptstyle{\ds}$ 
expriment deux familles de relations alg\'ebriques entre les valeurs
z\^eta multiples, qui engendrent conjecturalement toutes les relations. 
En analogie avec cette conjecture, les \'el\'ements de l'alg\`ebre de
double m\'elange elliptique $\scriptstyle{\ds_{ell}}$ sont des 
polyn\^omes de Lie ayant une propri\'et\'e dimorphique, appel\'ee
$\scriptstyle{\Delta}$-bialternalit\'e, qui d\'ecrit conjecturalement
(le dual de) l'ensemble des relations alg\'ebriques entre les
{\weeit valeurs z\^etas elliptiques multiples}, qui sont les coefficients
d'une certaine s\'erie g\'en\'eratrice elliptique (construite explicitement
dans [LMS]) reli\'ee \`a l'associateur d'Enriquez ([En1]).  Nous montrons
que l'un des r\'esultats majeurs de la th\'eorie des moules de J.~\'Ecalle
peut \^etre interpr\'et\'e comme l'existence d'un morphisme injectif
$\scriptstyle{\ds\rightarrow \ds_{ell}}$ d'alg\`ebres de Lie. 
Notre r\'esultat principal est la compatibilit\'e de ce morphisme avec
la section ``point base tangentiel'' 
$\scriptstyle{{\rm Lie}\,\pi_1(MTM)\rightarrow {\rm Lie}\,\pi_1(MEM)}$ 
construite par Hain and Matsumoto ([HM]), et avec 
la section $\scriptstyle{\grt\rightarrow \grt_{ell}}$ construite par
Enriquez qui envoie l'alg\`ebre de Lie $\scriptstyle{\grt}$ 
de Grothendieck-Teichm\"uller vers sa version elliptique
$\scriptstyle{\grt_{ell}}$.  Ces compatibilit\'es sont exprim\'ees par
la commutativit\'e du diagramme (A) (\`a l'exception de la fl\`eche
en pointill\'e, qui est conjecturale.)}}\par}
\vfill\eject
\noindent {\bf 1. Introduction}
\vskip .2cm
\noindent {\bf 1.1. Overview}
\vskip .1cm
The goal of this paper is to apply Ecalle's mould theory
to define an {\it elliptic double shuffle} Lie algebra $\ds_{ell}$ that 
turns out to parallel Enriquez' construction in [En1] of the elliptic 
Grothendieck-Teichm\"uller Lie algebra, and Hain and Matsumoto's construction  
of the fundamental Lie algebra of the category $MEM$ of mixed elliptic motives 
in [HM].  Both of those Lie algebras are equipped with canonical surjections 
to the corresponding genus zero Lie algebras, 
$$\cases{\grt_{ell}\rightarrow\!\!\!\!\rightarrow \grt\cr
{\rm Lie}\,\pi_1(MEM)\rightarrow\!\!\!\!\rightarrow{\rm Lie}\,\pi_1(MTM).}$$
\vskip .2cm
Here, $MTM$ is the category of mixed Tate motives over ${\Bbb Z}$, and
the notation ${\rm Lie}\,\pi_1(MTM)$ (resp.~${\rm Lie}\,\pi_1(MEM)$) denotes
the Lie algebra of the pro-unipotent radical of the fundamental
group of the Tannakian category $MTM$ (resp.~$MEM$) equipped with the
de Rham fiber functor (resp.~its lift to a fiber functor on $MEM$ via
composition with the natural surjection $MEM\rightarrow MTM$, cf. [HM,
\S 5].)

Each of the Lie algebras $\grt_{ell}$ and ${\rm Lie}\,\pi_1(MEM)$ is also
equipped with a natural section of the above surjection, corresponding, 
geometrically, to the tangential
base point at infinity on the moduli space of elliptic curves:
$$\cases{\gamma:\grt\hookrightarrow\grt_{ell}\cr
\gamma_t:{\rm Lie}\,\pi_1(MTM)\hookrightarrow{\rm Lie}\,\pi_1(MEM).}$$
Hain-Matsumoto determine a canonical Lie ideal of
$\frak{u}$ of ${\rm Lie}\,\pi_1(MEM)$, and Enriquez defines a canonical
Lie ideal $\r_{ell}$ of $\grt_{ell}$, such that the above sections
give semi-direct product structures
$$\cases{\grt_{ell}\simeq \r_{ell}\ {\Bbb o}\ \gamma(\grt)\cr
{\rm Lie}\,\pi_1(MEM)\simeq \frak{u}\ {\Bbb o}\ \gamma_t\bigl({\rm Lie}\,\pi_1(MTM)\bigr).}$$
\vskip .2cm
Let $C_i=ad(a)^{i-1}(b)$ for $i\ge 1$, and
let ${\rm Lie}[C]$ denote the Lie algebra ${\rm Lie}[C_1,C_2,\ldots]$.
It is an easy consequence of Lazard elimination that ${\rm Lie}[C]$ is 
a free Lie algebra on the generators $C_i$, and that
$${\rm Lie}[a,b]\simeq \Q a\oplus {\rm Lie}[C]$$
(see Appendix).  In other words, the elements in ${\rm Lie}[C]$ are
all the elements of ${\rm Lie}[a,b]$ having no linear term in $a$.
\vskip .3cm
\noindent {\bf Definition.} 
Let $\Der^0{\rm Lie}[a,b]$ denote the subspace of derivations $D\in 
\Der\,{\rm Lie}[a,b]$ that annihilate $[a,b]$ and 
such that $D(a)$ and $D(b)$ lie in ${\rm Lie}[C]$.

Hain-Matsumoto and Enriquez both give derivation 
representations of the elliptic spaces into $\Der^0{\rm Lie}[a,b]$, but
Enriquez's Lie morphism $\grt_{ell}\rightarrow \Der^0{\rm Lie}[a,b]$
is injective (by [T2], cf. below for more detail), 
whereas Hain-Matsumoto conjecture this result in the
motivic situation.  However, Hain-Matsumoto
compute the image of $\frak{u}$ in $\Der^0{\rm Lie}[a,b]$ and show that
it is equal to a certain explicitly determined Lie algebra $\frak{b}_3$
related to ${\rm SL}_2({\Bbb Z})$ (or to the Artin braid group $B_3$ on three
strands), namely the Lie algebra generated by derivations 
$\epsilon_{2i}$, $i\ge 0$ defined by $\epsilon_{2i}(a)=ad(a)^{2i}(b)$,
$\epsilon_{2i}([a,b])=0$\footnote{$^1$}{\wee This Lie algebra was 
introduced by Tsunogai in [T1, \S 3] (see also [P], [BS] and [Br4]
for some results on its interesting structure. The $\epsilon_{2i}$ also
play an important role in [CEE] and [En1].}, whereas Enriquez 
considers the same Lie algebra $\frak{b}_3$,
shows that it injects into $\r_{ell}$,
and conjectures that they are equal\footnote{$^2$}{\wee It is
really remarkable that these two papers were written totally independently
of one another.}. 

All these maps are compatible with the canonical injective morphism 
${\rm Lie}\,\pi_1(MTM)\rightarrow \grt$ whose existence was proven by 
Goncharov and Brown in two stages as follows.  Goncharov constructed a 
Hopf algebra ${\cal A}$ of motivic zeta values as motivic iterated integrals 
[G,\S 5], and identified it with a subalgebra of the Hopf algebra of framed 
mixed Tate motives [G, \S 8]; he showed that these motivic zeta 
values satisfy the associator relations. Brown [Br1] subsequently lifted
Goncharov's construction to an algebra ${\cal H}$ in which the motivic 
$\zeta^m(2)$ is non-zero, such that in fact ${\cal H}\simeq {\cal A}\otimes 
\Q[\zeta^m(2)]$.  He was able to compute the structure and the dimensions of
the graded parts of ${\cal H}$ and thus of ${\cal A}$, from which it follows 
that ${\cal A}$
is in fact equal to the full Hopf algebra of framed mixed Tate motives.  
In the dual situation, this means that the fundamental Lie
algebra of $MTM$ injects into the Lie algebra of associators, namely the
top arrow of the following commutative diagram:
$$\xymatrix{{\rm Lie}\,\pi_1(MTM)\ar[d]\ar[r]&\grt\ar[d]\\
{\rm Lie}\,\pi_1(MEM)\ar[r]&\grt_{ell}.}$$
\vskip .2cm
The elliptic double shuffle Lie algebra $\ds_{ell}$ that
we define in this article is conjecturally isomorphic to ${\rm Lie}\,\pi_1
(MEM)$ and $\grt_{ell}$.  We show that it shares with them the following 
properties: firstly, it comes equipped with an injective Lie algebra morphism 
$$\gamma_s:\ds\rightarrow \ds_{ell},$$ 
where $\ds$ is the {\it regularized} double shuffle Lie algebra 
defined in [R], where it is denoted $\dmr$ 
(``double m\'elange r\'egularis\'e'').

Secondly there is an injective derivation representation 
$$\ds_{ell}\hookrightarrow \Der^0{\rm Lie}[a,b].$$  
Unfortunately, we have not yet been able to find a good canonical Lie ideal in
$\ds_{ell}$ that would play the role of $\frak{u}$ and $\r_{ell}$, although
it is easy to show that there is an injection $\frak{b}_3\hookrightarrow
\ds_{ell}$ whose image conjecturally plays this role (cf. the end of
section 1.3).  Since 
$\frak{u}\rightarrow\frak{b}_3\hookrightarrow \ds_{ell}$, we do have
a Lie algebra injection,
$${\rm Lie}\,\pi_1(MEM)\hookrightarrow \ds_{ell},$$
but not the desired injection
$$\grt_{ell}\hookrightarrow\ds_{ell},$$
(the dotted arrow in the diagram in the abstract), which would follow
as a consequence of Enriquez' conjecture that $\r_{ell}=\frak{b}_3$.
It would have been nice to give a direct proof of the existence of a 
Lie algebra morphism $\grt_{ell}\rightarrow \ds_{ell}$ even without proving
Enriquez' conjecture, but we were not able to find one.  This result appears 
like an elliptic version of Furusho's injection $\grt\hookrightarrow \ds$ 
(cf. [F]), and may possibly necessitate some similar techniques.  

Our main result, however, is the commutation of the diagram given in the 
abstract, which does not actually require an injective map 
$\grt_{ell}\rightarrow \ds_{ell}$, but, given all the observations above,
comes down to the commutativity of the triangle diagram
$$\xymatrix{\grt\ar[rr]\ar[dr]&&\ds\ar[dl]\\
&\Der^0{\rm Lie}[a,b]&.}\eqno(1.1.1)$$
The morphisms from $\grt$ and $\ds$ to $\Der\,{\rm Lie}[a,b]$ factor
through the respective elliptic Lie algebras (cf. the diagram in the abstract).
Note that the morphisms in (1.1.1) must not be confused with the familiar
Ihara-type morphism
$\grt\rightarrow \Der\,{\rm Lie}[x,y]$ via
$y\mapsto [\psi(-x-y,y),y]$ and
$x+y\mapsto 0$, and the analogous map for $\ds$ investigated in [S2].
The relation between the two is based on the fact that ${\rm Lie}[x,y]$
is identified with the Lie algebra of the fundamental group of the 
thrice-punctured sphere, whereas ${\rm Lie}[a,b]$ is identified with the
Lie algebra of the once-punctured torus.  
The natural Lie morphism ${\rm Lie}[x,y]\rightarrow {\rm Lie}[a,b]$,
reflecting the underlying topology, is given by
$$x\mapsto t_{01}, y\mapsto t_{02},$$
where we write $Ber_x=ad(x)/\bigl(exp(ad(x))-1\bigr)$
for any $x\in{\rm Lie}[a,b]$, and set 
$$t_{01}=Ber_{b}(-a),\ \ t_{02}=Ber_{-b}(a).$$
We show that certain derivations of ${\rm Lie}[x,y]$,
transported to the free Lie subalgebra ${\rm Lie}[t_{01},t_{02}]
\subset {\rm Lie}[a,b]$ 
have a unique extension to derivations of all of ${\rm Lie}[a,b]$, and that 
in particular this is the case for the derivations in the image of 
$\grt$ and $\ds$ (cf. section 2). This gives a direct
interpretation of the two maps to derivations in the diagram (1.1.1) whose
commutativity we prove.

The existence of the injection 
$\ds\rightarrow\ds_{ell}$ arose from an elliptic reinterpretation of 
a major theorem by Ecalle in mould theory. This reading of Ecalle's work 
and interpretation of some of his important results constitute one of the 
main goals of this paper in themselves.   Indeed, it appears that Ecalle's 
seminal work in mould and multizeta theory has been largely ignored by the 
multiple zeta community\footnote{$^3$}{\wee According to
the author's discussion with several colleagues, this appears to be at least
partly due to a reluctance to accept Ecalle's language, because, at least
according to some, it uses a system of words with varying vowels, rather 
than the more standard single letters, for the basic objects. This seems surprising,
as it is unclear why calling a derivation arit(f), say, rather than 
D$_f$ should pose such a problem.  Possibly we enter here into the domain
of psychology.  A second, more serious obstacle is the lack of
proofs in Ecalle's work, and the incredible profusion of statements,
which makes it difficult to pick out exactly what is needed to establish
a specific result. The author has attempted to solve this problem, at least
partially, in the basic text [S] which gives an introduction with complete 
proofs to the portion of Ecalle's work most directly related to current 
problems in double shuffle algebra.}.
This minimalist way of phrasing the main result shows that it could actually 
be stated and proved without even defining an elliptic double shuffle Lie 
algebra.  However, this object is important in its own right, principally for 
the following reason. Recall that the usual double shuffle Lie algebra $\ds$ 
expresses the double shuffle relations satisfied by the multiple zeta values, 
in the following sense.  Let ${\cal FZ}$, the {\it formal multizeta
algebra}, be the graded dual of the universal enveloping algebra of $\ds$;
it is generated by formal symbols satisfying only the double shuffle 
relations. Since motivic and real multizeta values are known to satisfy
them (see for example [So]), ${\cal FZ}$ surjects onto the algebras of motivic and real
multizeta values. These surjections are conjectured to be isomorphisms,
i.e., it is conjectured that the double shuffle relations generate all
algebraic relations between motivic, resp.~real multizeta values (with the
first of these problems being undoubtedly much more tractable than the
second, for reasons of transcendence).

The role played by the double shuffle algebra with respect to ordinary
multizeta values is analogous to the role played by the elliptic double
shuffle algebra defined in this article with respect to the elliptic
mzv's defined in [LMS]. There, we define an {\it elliptic generating series}
in the completed Lie algebra ${\rm Lie}[a,b]$, whose coefficients,
called elliptic mzv's or emzv's, are related to the iterated integrals that 
form the coefficients of Enriquez' monodromic elliptic associator, and we give an 
explicit ``dimorphic'' or ``double shuffle'' type symmetry of this 
generating series which is exactly the defining property of $\ds_{ell}$.  
Indeed, letting ${\cal E}$ denote the graded Hopf algebra generated by 
the emzv's, we show in [LMS] that the vector space 
$\frak{ne}={\cal E}_{>0}/\bigl({\cal E}_{>0}\bigr)^2$  
is isomorphic to a semi-direct product $\frak{b}_3{\Bbb o}\nz^\vee$, where
$\nz$ is the space of ``new multizeta values'' obtained by quotienting
the algebra of multizeta values by $\zeta(2)$ and products.
Under the standard conjecture from multizeta theory $\nz^\vee\simeq\grt$,
as well as Enriquez' conjecture $\frak{r}_{ell}\simeq \frak{b}_3$,
this implies that $\frak{ne}\simeq \grt_{ell}$.  If $\grt_{ell}\simeq \ds_{ell}$,
as we believe, this would mean that the elliptic double shuffle property 
determines all algebraic relations between the emzv's.
This topic, which reflects the geometric 
aspects of the elliptic double shuffle relations introduced in this paper, 
is explored in detail in [LMS].

The content of the present paper has some relation with the recent preprint
[Br3] as well as the earlier, closely related online lecture notes [Br2]. 
In particular Brown gives the existence of rational-function moulds 
satisfying the double shuffle relations, which is an immediate consequence 
of an important theorem of Ecalle that appears in all of his articles 
concerning ARI/GARI and multiple zeta
values (cf.~Theorem 1.3.2 below), although Brown introduces a completely 
different construction (vines and grapes). Brown also mentions in passing
(cf. (3.7) of [Br3]) the result of the useful extension Lemma 2.1.2 below, 
however without proof. In [Br2] (conjecture 3) and [Br3]
(following Prop. 4.6), Brown asks the question of whether
$\u^{geom}\simeq \pls$.  The answer to this question is no; indeed
all elements of $\grt$ with no depth 1 part furnish elements of $\pls$ not 
lying in $\frak{u}$ via Enriquez' section, as explained in the Corollary 
following Theorem 1.3.3.

\vskip .2cm
\noindent {\bf Acknowledgements.} The work on this paper benefited from
discussions with B. Enriquez and P. Lochak, both of whom listened patiently
and provided some crucial elements of proof. J. Ecalle repeatedly gave of his 
time to help understand some of his results. R. Hain also shed some light on 
details arising from his motivic work. I thank them all warmly.  

\vskip .5cm
\noindent {\bf 1.2. The elliptic Grothendieck-Teichm\"uller Lie algebra }
\vskip .2cm
In this section we recall the definition of the elliptic 
Grothendieck-Teichm\"uller Lie algebra $\grt_{ell}$ defined in [En1], 
along with some of its main properties.
Recall that the genus 1 braid Lie algebra on $n$ strands, $\t_{1,n}$,
is generated by elements
$x_1^+,\ldots,x_n^+$ and $x_1^-,\ldots,x_n^-$ subject to relations
$$x_1^++\cdots+x_n^+=x_1^-+\cdots+x_n^-=0,\ \ 
[x_i^+,x_j^+]=[x_i^-,x_j^-]=0\ \ {\rm if\ i\ne j},$$
$$[x_i^+,x_j^-]=[x_j^+,x_i^-]\ \ {\rm for\ } i\ne j,\ \ 
[x_i^+,[x_j^+,x_k^-]]=[x_i^-,[x_j^+,x_k^-]]=0\ \ {\rm for}\ i,j,k\ 
{\rm distinct}.$$
We write $t_{ij}=[x_i^+,x_j^-]$.
It is the Lie algebra of the unipotent completion of the topological fundamental
group of the configuration space of $n$ ordered marked points on the torus
(cf.~[CEE,\S 2.2] for details).
The Lie algebra $\t_{1,2}$ is isomorphic to ${\rm Lie}[a,b]$, the free 
Lie algebra on two generators\footnote{$^4$}{\wee With 
respect to the notation of [En1] we have ${\rm Lie}[a,b]={\t}_{1,2}$, 
$a=y_1=x_1^-$, $b=x_1=x_1^+$ (Enriquez uses both notations).}
$a$ and $b$.
Throughout this article, we write ${\rm Lie}[a,b]$
for the completed Lie algebra, i.e., it contains infinite Lie series and not
just polynomials.  Thus an element $\alpha\in \t_{1,2}\simeq {\rm Lie}[a,b]$
is a Lie series $\alpha(a,b)$ in two free variables.

\vskip .3cm
\noindent {\bf Definition.} The {\it elliptic Grothendieck-Teichm\"uller Lie
algebra} $\grt_{ell}$ is the set of triples 
$(\psi,\alpha_+,\alpha_-)$
with $\psi\in \grt$, $\alpha_+,\alpha_-\in \t_{1,2}$, such that setting
$$\cases{\Psi(x_1^{\pm})=\alpha_{\pm}(x_1^{\pm},x_1^{\mp})
+[x_1^{\pm},\psi(t_{12},t_{23})],\cr
\cr
\Psi(x_2^{\pm})=\alpha_{\pm}(x_2^{\pm},x_2^{\mp})
+[x_2^{\pm},\psi(t_{12},t_{13})],\cr
\cr
\Psi(x_3^{\pm})=\alpha_{\pm}(x_3^{\pm},x_3^{\mp})}\eqno(1.2.1)$$
yields a derivation of ${\frak{t}}_{1,3}$. The space $\grt_{ell}$ is made
into a Lie algebra by bracketing derivations; in other words,
writing $D_{\alpha_{\pm}}$ for the derivation of $\t_{1,2}\simeq
{\rm Lie}[a,b]$ which takes
$a\mapsto \alpha_+(a,b)$ and $b\mapsto \alpha_-(a,b)$, we have
$$\langle (\psi,\alpha_+,\alpha_-),(\phi,\beta_+,\beta_-)\rangle=
\Bigl(\{\psi,\phi\},D_{\alpha_{\pm}}(\beta_+)-D_{\beta_{\pm}}(\alpha_+),
D_{\alpha_{\pm}}(\beta_-)-D_{\beta_{\pm}}(\alpha_-)\Bigr),$$
where $\{\psi,\phi\}$ is the Poisson (or Ihara) bracket on $\grt$.
Finally, we assume that the coefficient of $a$ in both 
$\alpha_+$ and $\alpha_-$ is equal to $0$.  
\vskip .3cm
\noindent {\bf Remark.} The last assumption is not contained
in Enriquez' original definition. In particular he allows the 
element $(0,0,a)$, corresponding to the derivation $e(a)=0$,
$e(b)=a$, which together with $\epsilon_0(a)=b$, $\epsilon_0(b)=0$
generate a copy of $\frak{sl}_2$ in $\grt_{ell}$.
Because of this, Enriquez' version of $\grt_{ell}$ is not pronilpotent, and
is thus strictly larger than the ${\rm Lie}\,\pi_1(MEM)$ studied in
[HM], which is the Lie algebra of the prounipotent radical of the fundamental
group of $MEM$. Thus, isomorphism can only be conjectured if the extra 
element is removed, motivating our slight 
alteration of his definition.  We nonetheless write $\grt_{ell}$ for the 
modified version; the results of Enriquez on elements of $\grt_{ell}$ that 
we cite adapt directly with no changes.
\vskip .2cm

We summarize Enriquez' important results concerning $\grt_{ell}$ in the
following theorem.
\vskip .3cm
\noindent {\bf Theorem 1.2.1.} (cf.~[En1]) {\it For all $(\psi,\alpha_+,\alpha_-)
\in \grt_{ell}$, the derivation $D_{\alpha_{\pm}}$ of $\t_{1,2}$ annihilates 
the element $t_{12}=[a,b]$.
But for each $\psi\in\grt$, there exists one and only one triple 
$(\psi,\alpha_+,\alpha_-)\in \grt_{ell}$ such that
$D_{\alpha_{\pm}}$ restricts to the Lie subalgebra ${\rm Lie}[t_{01},t_{12}]$
as follows:
$$\cases{t_{01}\mapsto [\psi(t_{01},t_{12}),t_{01}],\cr
t_{02}\mapsto [\psi(t_{02},t_{12}),t_{02}],\cr
t_{12}\mapsto 0.}\eqno(1.2.2)$$
The map $\gamma:\grt\rightarrow \grt_{ell}$ mapping $\psi$ to
this triple is a Lie algebra morphism that is a section of the 
canonical surjection $\grt_{ell}\rightarrow \grt$.  The Lie algebra $\grt_{ell}$ thus has a 
semi-direct product structure
$$\grt_{ell}=\r_{ell}{\Bbb o}\gamma(\grt).\eqno(1.2.3)$$}
These results of Enriquez show that $\grt_{ell}$ is generated by elements 
belonging to two particular subspaces.  The first is the subspace $\r_{ell}$ of 
triples $(\psi,\alpha_+,\alpha_-)$ with $\psi=0$, which forms a Lie ideal inside
$\grt_{ell}$.  The quotient $\grt_{ell}/\r_{ell}$ is canonically isomorphic
to $\grt$, the surjection being nothing other than the morphism 
forgetting $\alpha_+$ and $\alpha_-$.  The second subspace, the image
of the section $\grt\hookrightarrow\grt_{ell}$, is the space of 
triples that restrict on the free Lie subalgebra ${\rm Lie}[t_{01},t_{02}]$ to 
Ihara-type derivations (1.2.2).
For any triple $(\psi,\alpha_+,\alpha_-)$ of the second type,
i.e., in the -- but only
and uniquely for those, not for general elements of $\grt_{ell}$ --
we let $D_\psi=D_{\alpha_{\pm}}$, and write $\tilde D_\psi$ for
the the restriction of $D_\psi$ to ${\rm Lie}[t_{02},t_{12}]$ given by (1.2.2). 
\vskip .2cm
\noindent {\bf Remark.} This is actually a rephrasing of part of Enriquez' results.  In fact, he 
gives the derivation $D_\psi$ by explicitly displaying its value on $t_{01}$ 
(as in (1.2.2) and on $b$.  Since $D_\psi(t_{12})=0$, the restriction of 
$D_\psi$ to ${\rm Lie}[t_{01},t_{02}]$ is the well-known Ihara derivation 
associated to $\psi\in \grt$, and therefore the value on $t_{02}$ must be 
as in (1.2.2). The fact that $D_\psi$ is the only extension of (1.2.2) to a 
derivation on all of ${\rm Lie}[a,b]$ follows from our extension Lemma 2.1.2 
below.  This characterization of $D_\psi$ is sufficient for our purposes
in this article; we do not actually use the explicit expression of
$D_\psi(b)$, but it is necessary for Enriquez' work on elliptic
associators.
\vskip .2cm
The map 
$$\eqalign{\grt_{ell}&\rightarrow \Der^0{\rm Lie}[a,b]\cr
(\psi,\alpha_+,\alpha_-)&\mapsto D_{\alpha_{\pm}}}$$
is injective; in other words, knowing the pair $(\alpha_+,\alpha_-)$ allows 
us to uniquely recover $\psi$.  This result follows from [T2, Theorem 1.17]
(building on previous work in [NTU]), which states that removing the third
braid strand yields an injection ${\cal D}_1^{(2)}\hookrightarrow 
{\cal D}_1^{(1)}$, where ${\cal D}_1^{(1)}\simeq \Der^0{\rm Lie}[a,b]$ and
${\cal D}_1^{(2)}$ is a space of special derivations of ${\cal L}_1^{(2)}
\simeq \t_{1,3}$ which contains (and is conjecturally equal to) $\grt_{ell}$.

Furthermore, by Lemma 2.1.1 below, there is an injective linear map 
$$\eqalign{\Der^0{\rm Lie}[a,b]&\rightarrow {\rm Lie}[a,b]\cr
D&\mapsto D(a),}\eqno(1.2.4)$$
which is a Lie algebra bijection onto its image when that image (equal
to the subspace ${\rm Lie}^{push}[a,b]$ of {\it push-invariant} elements of ${\rm Lie}[a,b]$, 
cf. section 2) is equipped with the corresponding bracket. In particular
this shows that in the triple $(\psi,\alpha_+,\alpha_-)$, the
element $\alpha_+$ determines $\alpha_-$, and thus also $\psi$.
We write $\gamma_+:\grt\hookrightarrow {\rm Lie}[a,b]$ for the
map sending $\psi\mapsto \alpha_+$.  By the above arguments,
$\gamma_+$ determines $\gamma$ and vice versa. 

The desired triangle diagram (1.1.1) is thus equivalent to
$$\xymatrix{\grt\ar[dr]^{\gamma_+}\ar@{^(->}[rr]&& \ds\ar[dl]\\
&{\rm Lie}^{push}[a,b],&}\eqno(1.2.5)$$
by composing it with the map (1.2.4).  Our main result, Theorem 1.3.1 below,
is the explicit version of the commutation of the diagram (1.2.5).
\vskip .5cm
\noindent {\bf 1.3. Mould theory, elliptic double shuffle and the main theorem}
\vskip .3cm
In this section we explain how we use Ecalle's mould theory -- particularly
adapted to the study of dimorphic (or ``double shuffle'') structures --
to construct
the {\it elliptic double shuffle} Lie algebra $\ds_{ell}$, which like
$\grt_{ell}$ is a subspace of the push-invariant elements of ${\rm Lie}[a,b]$,
and how we reinterpret one of Ecalle's major theorems and combine it with
some results from Baumard's Ph.D. thesis ([B]), to define the injective 
Lie morphism $\ds\rightarrow \ds_{ell}$.

We assume some familiarity with moulds in this section; however all the 
necessary notation and definitions starting with that of a mould are recalled 
in the appendix at the end of the paper.  We use the notation $ARI$ to denote 
the vector space of moulds with constant term 0, and write $ARI_{lu}$ for $ARI$ 
equipped with the $lu$-bracket and $ARI_{ari}$ for $ARI$ equipped with the
$ari$-bracket (the usual $ARI$ according to Ecalle's notation).  Simi\-larly,
we write $GARI$ for the set of moulds with constant term 1 and write 
$GARI_{mu}$ and $GARI_{gari}$ for the groups obtained by equipping $GARI$ with 
the $mu$ and $gari$ multiplication laws.  In section 3 we will introduce a 
third Lie bracket on $ARI$, the $Dari$-bracket, and employ the notation 
$ARI_{Dari}$, as well as the corresponding group $GARI_{Dgari}$ with 
multiplication law $Dgari$.  
\vskip .2cm
We define the following operators on moulds:
$$\cases{dar(P)(u_1,\ldots,u_r)=u_1\cdots u_r\,P(u_1,\ldots,u_r)\cr
dur(P)(u_1,\ldots,u_r)=(u_1+\cdots+u_r)\,P(u_1,\ldots,u_r)\cr
\Delta(P)(u_1,\ldots,u_r)=u_1\cdots u_r(u_1+\cdots+u_r)\,P(u_1,\ldots,u_r)\cr
ad(Q)\cdot P=[Q,P]\ \ {\rm for\ all\ }Q\in ARI.}\eqno(1.3.1)$$
\vskip .2cm
We take $dar(P)(\emptyset)= dur(P)(\emptyset)=\Delta(P)(\emptyset)=
P(\emptyset)$. The operators $dur$ and $ad(Q)$
are derivations of the Lie algebra $ARI_{lu}$, whereas $dar$ is an 
automorphism of $ARI_{lu}$.
We will also make use of the inverse operators $dur^{-1}$ 
(resp. $dur^{-1}$ and $\Delta^{-1}$) defined by dividing a mould in depth 
$r$ by $(u_1+\cdots+u_r)$ (resp. by $(u_1+\cdots+u_r)$ and $(u_1+\cdots+u_r)u_1
\cdots u_r$).  
\vskip .2cm

If $p\in {\rm Lie}[a,b]$, then we have
$$\cases{ma\bigl([p,a]\bigr)=dur\bigl(ma(p)\bigr)\cr
ma\bigl(p(a,[b,a])\bigr)=dar\bigl(ma(p)\bigr)\cr
ma\bigl([p(a,[b,a]),a]\bigr)=\Delta\bigl(ma(p)\bigr).}\eqno(1.3.2)$$
A proof of the first equality can be found in [R, Proposition 4.2.1.1] 
or [S, Lemma 3.3.1]. The second is obvious from the definition of $ma$
(cf. Appendix), since substituting $[b,a]$ for
$b$ in $C_k$ yields $-C_{k+1}$ so making the substitution in
a monomial $C_{k_1}\cdots C_{k_r}$ yields $(-1)^r C_{k_1+1}\cdots C_{k_r+1}$,
and we have
$$ma\bigl((-1)^r C_{k_1+1}\cdots C_{k_r+1}\bigr)=(-1)^r(-1)^{k_1+\cdots+k_r}
u_1^{k_1}\cdots u_r^{k_r}=u_1\cdots u_r\,ma\bigl(C_{k_1}\cdots C_{k_r}\bigr).$$
The third equality of (1.3.2) follows from the first two.
\vskip .2cm
We now recall the definition of the key mould $pal$ that lies at the heart
of much of Ecalle's theory of moulds.  Following [E2], we
start by introducing an auxiliary mould $dupal\in ARI$, given by
the simple explicit expression
$$dupal(u_1,\ldots,u_r)={{B_r}\over{r!}}
{{1}\over{u_1\cdots u_r}}\biggl(\sum_{j=0}^{r-1}
(-1)^j\Bigl({{r-1}\atop{j}}\Bigr)u_{j+1}\biggr).\eqno(1.3.3)$$
The mould $pal$ is then defined by setting $pal(\emptyset)=1$ and using
the equality
$$dur(pal)=pal\,dupal,\eqno(1.3.4)$$
which gives a recursive definition for $pal$ depth by depth starting with 
$pal(\emptyset)=1$, since to determine the left-hand side $dur(pal)$ 
in depth $r$ only requires knowing $pal$
up to depth $r-1$ on the right-hand side.

Since $pal(\emptyset)=1$, we have $pal\in GARI$.  We write $invpal$ for its 
inverse $inv_{gari}(pal)$ in the group $GARI_{gari}$.  Since
$GARI_{gari}$ is the exponential of the Lie algebra $ARI_{ari}$, it has an
adjoint action on $ARI_{ari}$; we write $Ad_{ari}(P)$ for the adjoint operator
on $ARI_{ari}$ associated to a mould $P\in GARI_{gari}$. 

At this point we are already equipped to baldly state our
main theorem linking Ecalle's theory of moulds to Enriquez' section
$\gamma:\grt\rightarrow \grt_{ell}$, or rather to the modified version
$\gamma_+$ introduced above that maps $\psi$ to the associated element
$\alpha_+$ in Enriquez' triple $(\psi,\alpha_+,\alpha_-)$.
\vskip .3cm
\noindent {\bf Theorem 1.3.1.} {\it Let $\psi\in\grt$ and set $f(x,y)=\psi(x,-y)$.  We have the following equality of moulds:
$$\Delta\bigl(Ad_{ari}(invpal)\cdot ma(f)\bigr)=ma\bigl(\gamma_+(\psi)\bigr).\eqno(1.3.5)$$}

In order to place this theorem in context and explain its power in
terms of helping to define an elliptic double shuffle Lie algebra that in 
turn will shed light on the dimorphic (``double-shuffle'') properties of
elliptic multiple zeta values, we first give some results from the literature,
starting with Ecalle's main theorem, with which he first revealed the
surprising role of the adjoint operator $Ad_{ari}(pal)$ and its
inverse $Ad_{ari}(pal)^{-1}=Ad_{ari}(invpal)$.

Recall from the appendix that in terms of moulds, $\ds$ is isomorphic
to the Lie subalgebra of $ARI_{ari}$ of polynomial-valued moulds that are
even in depth 1, and are alternal with swap that is alternil up to addition of 
a constant mould.  The notation we use for this in mould
language is a bit heavy, but has the advantage of concision and total 
precision in that the various symbols attached to $ARI$ carry all of the 
information about the moulds in the subspace under consideration: we have
the isomorphism
$$ma:\ds\buildrel{ma}\over\rightarrow ARI_{ari}^{pol,\underline{al}*\underline{il}},$$
where {\it pol} indicates polynomial moulds, the underlining is Ecalle's 
notation for moulds that are even in depth 1, and the usual notation 
$al/il$ for an alternal mould with alternil swap is weakened to 
$al*il$ when the swap is only alternil up to addition of a constant mould.

Similarly, the notation $ARI_{ari}^{\underline{al}*\underline{al}}$ refers to
the sub§space of moulds in $ARI_{ari}$ that are even in depth 1 and 
alternal with swap that is alternal up to addition of a constant mould 
(or ``bialternal'').  When we consider the subspace of these moulds that 
are also polynomial-valued, $ARI^{pol,\underline{al}*\underline{al}}$,
we obtain the (image under $ma$ of the) ``linearized double shuffle''
space $\ls$ studied for example in [Br3]. But the full non-polynomial space 
is of course hugely larger.  One of Ecalle's most remarkable discoveries is 
that the mould $pal$ provides an isomorphism between the two types of dimorphy,
as per the following theorem.
\vskip .3cm
\noindent {\bf Theorem 1.3.2.} (cf.~[E]\footnote{$^5$}{\wee This result is stated and used constantly in [E], as well as many other analogous results concerning
other symmetries.  But the proof is not given.  Ecalle was kind enough to
send us a sketch of the proof in a personal letter, relying on the fundamental
identity (2.62) of [E], itself not proven there.  Full details of the
reconstructed proof can be found in [S], with (2.62) proved in Theorem 2.8.1
and Theorem 1.3.2 above proved in \S 4.6.})
 {\it The adjoint map $Ad_{ari}(invpal)$ 
induces a Lie isomorphism of Lie subalgebras of $ARI_{ari}$:
$$Ad_{ari}(invpal):ARI^{\underline{al}*\underline{il}}_{ari}\buildrel{\sim}\over\rightarrow 
ARI^{\underline{al}*\underline{al}}_{ari},\eqno(1.3.6)$$
and if $F\in ARI^{\underline{al}*\underline{il}}$ and $C$ is the constant
mould such that $swap(F+C)$ is alternil, then 
$swap\bigl(Ad_{ari}(invpal)(F)\bigr)+C$ is alternal, i.e., the
constant corrections for $F$ and $Ad_{ari}(invpal)\cdot F$ are the same.  
In particular if $C=0$, i.e., if $F$ is $\underline{al}/\underline{il}$, then 
$Ad_{ari}(invpal)(F)$ lies in $\underline{al}/\underline{al}$.  }
\vskip .3cm
One important point to note in the result of Theorem 1.3.2 is that
the operator $Ad_{ari}(invpal)$ does not respect polynomiality of moulds.
Indeed, applying $Ad_{ari}(pal)$ to bialternal polynomial moulds produces
quite complicated denominators with many factors.  However, in his
doctoral thesis S. Baumard was able to show that conversely, when 
applying $Ad_{ari}(invpal)$ to moulds $ma(f)$ for $f\in \ds$, i.e.,~to
moulds in $ARI^{pol,\underline{al}*\underline{il}}$, the
denominators remain controlled.  Indeed, let $ARI^\Delta$
denote the space of moulds $P\in ARI$ such that $\Delta(P)\in ARI^{pol}$,
i.e., the space of rational-function valued moulds whose denominator 
is ``at worst'' $u_1\cdots u_r(u_1+\cdots+u_r)$ in depth $r$. 
\vskip .3cm
\noindent {\bf Theorem 1.3.3.} [B, Thms. 3.3, 4.35] {\it The space $ARI^\Delta$
forms a Lie algebra under the $ari$-bracket, and we have an injective Lie
algebra morphism
$$Ad_{ari}(invpal):ARI^{pol,\underline{al}*\underline{il}}_{ari}\hookrightarrow ARI_{ari}^\Delta.\eqno(1.3.7)$$}
Recall that $\pls$ (``polar linearized double shuffle'')
is the notation used by F. Brown for the space $ARI^{\Delta,\underline{al}/
\underline{al}}$ and $\frak{u}$ for the Lie subalgebra of $ARI$ generated
by $B_{-2}$ and $B_{2i}$ for $i\ge 1$, where $B_i$ denotes the mould
concentrated in depth 1 defined by $B_i(u_1)=u_1^i$.
As a corollary of Theorems 1.3.1, 1.3.2 and 1.3.3, we give a negative
answer to the question posed by Brown ([Br2], 
conjecture 3 and [Br3], following Prop. 4.6) as to whether $\pls$
and $\frak{u}$ are equal.
\vskip .3cm
\noindent {\bf Corollary.} {\it Let $\psi\in\grt$ be an element of
$\grt$ having no depth 1 part.  Then
$$\Delta^{-1}\bigl(ma(\gamma_+(\psi)\bigr)\in ARI^{\Delta,\underline{al}/\underline{al}}
=\pls$$
but
$$\Delta^{-1}\bigl(ma(\gamma_+(\psi))\bigr)\notin \frak{u}.$$}
\noindent Proof. Since by Furusho's theorem, $\psi(x,y)\mapsto
f(x,y)=\psi(x,-y)$ maps $\grt\hookrightarrow\ds$, we have 
$ma(f)\in ARI^{pol,\underline{al}*\underline{il}}$ for every $\psi\in\grt$.
In particular, if $\psi$ has no depth 1 part, then
$ma(f)\in ARI^{pol,\underline{al}/\underline{il}}$; thus by Theorem 1.3.2,
$Ad_{ari}(invpal)\cdot ma(f)\in ARI^{\underline{al}/\underline{al}}$,
and by Theorem 1.3.3, it also lies in $ARI^{\Delta}$; thus it lies in
$ARI^{\Delta,\underline{al}/\underline{al}}=\pls$. By Theorem 1.3.1,
$Ad_{ari}(invpal)\cdot ma(f)$ is equal to $\Delta^{-1}\bigl(ma(\gamma_+(\psi)\bigr)$ where $\gamma_+$ denotes Enriquez' section $\grt\rightarrow \grt_{ell}$,
associating to $\psi\in\grt$ the element $\alpha_+$ from the triple
$(\psi,\alpha_+,\alpha_-)$. But Enriquez shows that
$\grt_{ell}$ is a semi-direct product $\gamma_+(\grt){\Bbb o}\r_{ell}$
and that $\Delta(\frak{u})\subset ma(\r_{ell})$.  Thus
$ma\bigl(\gamma_+(\grt)\bigr)\cap \Delta(\frak{u})=\{0\}$.
\hfill{$\diamondsuit$}
\vskip .3cm
For the rest of this article we will use the notation:
$$\cases{f=\psi(x,-y)\cr
F=ma(f)\cr
A=Ad_{ari}(invpal)\cdot F\cr
M=\Delta(A).}\eqno(1.3.8)$$ 
\vskip .3cm
\noindent {\bf Corollary 1.3.4.} {\it Let $f\in \ds$ and let $F=ma(f)$, so
$F\in ARI^{pol,\underline{al}*\underline{il}}$. Then the mould 
$M=\Delta\bigl(Ad_{ari}(invpal)\cdot F\bigr)$ is alternal, push-invariant
and polynomial-valued.}
\vskip .2cm
\noindent Proof. Let $A=Ad_{ari}(invpal)\cdot F$. Then $A\in
ARI^{\underline{al}*\underline{al}}$ by Theorem 1.3.2, so $A$ is alternal,
and furthermore $A$ is push-invariant because all moulds in 
$ARI^{\underline{al}*\underline{al}}$ are push-invariant (see [E2] or 
[S, Lemma 2.5.5]). Thus $M=\Delta(A)$ is also alternal
and push-invariant since $\Delta$ preserves these properties.  The fact that 
$M$ is polynomial-valued follows from Theorem 1.3.3.\hfill{$\diamondsuit$}
\vskip .3cm
\noindent {\bf Definition.} A mould $P$ is said to be {\it $\Delta$-bialternal} 
if $\Delta^{-1}(P)$ is bialternal, i.e., $P\in \Delta(ARI_{ari}^{al*al})$.
The {\it elliptic double shuffle Lie algebra} $\ds_{ell}\subset 
{\rm Lie}[a,b]$ is the set of Lie polynomials which map under $ma$ to 
polynomial-valued $\Delta$-bialternal moulds that are even in depth 1, i.e.,
$$\ds_{ell}=ma^{-1}\Bigl(\Delta\bigl(ARI^{\Delta,\underline{al}*\underline{al}}_{ari}\bigr)\Bigr).\eqno(1.3.9)$$

Taken together, Theorems 1.3.2 and 1.3.3 show that the image of 
$ma(\ds)=ARI_{ari}^{pol,\underline{al}*\underline{il}}$ under $Ad_{ari}(invpal)$
lies in $ARI_{ari}^{\Delta,\underline{al}*\underline{al}}$, so the
image under $\Delta\circ Ad_{ari}(invpal)$ lies in the space of
polynomial-valued $\Delta$-bialternal moulds that are also even in depth 1
(since it is easy to see that $Ad_{ari}(invpal)$ preserves the lowest-depth
part of a mould). Thus we can define $\gamma_s$ to be the polynomial
avatar of $\Delta\circ Ad_{ari}(invpal)$, i.e., $\gamma_s$ is defined
by the commutation of the diagram
$$\xymatrix{\ds\ar[r]^{ma\ \ \ \ }\ar[d]_{\gamma_s}&ARI_{ari}^{pol,\underline{al}*\underline{il}}\ar[d]^{\Delta\circ Ad_{ari}(invpal)}\\
\ds_{ell}\ar[r]^{ma\ \ \ \ \ \ }&\Delta\bigl(ARI_{ari}^{\Delta,\underline{al}*\underline{al}}\bigr).}\eqno(1.3.10)$$
Thus for $f\in \ds$ we have 
$$ma\bigl(\gamma_s(f)\bigr)=\Delta\bigl(Ad_{ari}(invpal)\cdot ma(f)\bigr).$$
This reduces the statement of the main Theorem 1.3.1 above to the 
equality $$\gamma_s(f)=\gamma_+(\psi),$$
i.e., to the commutation of the diagram
$$\xymatrix{\grt\ar[rr]\ar[dr]_{\gamma_+}&&\ds\ar[dl]^{\gamma_s}\\
&{\rm Lie}^{push}[a,b],&}$$
which is the precise version of the desired diagram (1.2.5). 

As a final observation, we note that the definition of $\ds_{ell}$ 
makes the injective Lie algebra morphism $\frak{b}_3\hookrightarrow \ds_{ell}$ 
mentioned at the beginning of the introduction obvious.  Indeed, identifying
$\frak{b}_3$ with its image in ${\rm Lie}^{push}[a,b]$ under the map
(1.2.4), it is generated by the polynomials $\epsilon_{2i}(a)=ad(a)^{2i}(b)
=C_{2i+1}$ for $i\ge 0$, which map under $ma$ to the moulds $B_{2i}$ 
concentrated in depth 1 
and given by $B_{2i}(u_1)=u_1^{2i}$ (Ecalle denotes these moulds by 
$ekma_{2i}$ at least for $i\ge 1$; note however that $B_0$ and
$\Delta^{-1}(B_0)=B_{-2}$ are essential in the elliptic situation).  
To show that these moulds lie in $\ds_{ell}$, we need only note that
the moulds $\Delta^{-1}(B_{2i})=B_{2i-2}$ are even in depth 1, and 
trivially bialternal since this condition is empty in depth 1.
\vfill\eject
\noindent {\bf 2. Proof of the main theorem}
\vskip .3cm
For the proof of the main theorem, we first recall in 2.1 a 
few well-established facts about non-commutative polynomials, moulds and 
derivations, and give the key lemma about extending derivations on the Lie 
subalgebra ${\rm Lie}[t_{01},t_{02}]$ to all of ${\rm Lie}[a,b]$.  Once these 
ingredients are in place, the proof of the main theorem, given in 2.2, is a 
simple consequence of one important proposition, whose proof, contained in 
section 3, necessitates some developments in mould theory. In fact, the 
present section could be written entirely in terms of polynomials in $a$ 
and $b$ without any reference to moulds.  We only use moulds in the proof of 
Lemma 2.1.1, but merely as a convenience, as even this result could be stated 
and proved in terms of polynomials. Indeed this has already been done 
(cf. [S2]), but the proof given here using moulds is actually more 
elegant and simple.

\vskip .3cm
\noindent {\bf 2.1. The push-invariance and extension lemmas}
\vskip .3cm
\noindent {\bf Definition.} For $p\in {\rm Lie}[a,b]$, 
write $p=p_aa+p_bb$ and set 
$$p'=\sum_{i\ge 0}{{(-1)^{i-1}}\over{i!}} a^ib\partial_a^i(p_a)\eqno(2.1.1)$$ 
where $\partial_a(a)=1$, $\partial_a(b)=0$.  We call $p'$ the {\it partner} of 
$p$.  If $P\in ARI$ then we define $P'$ to be the {\it mould partner} of $P$,
given by the formula
$$P'(u_1,\ldots,u_r)={{1}\over{u_1+\cdots+u_r}}\Bigl(P(u_2,\ldots, u_{r-1},-u_1-\cdots-u_{r-1})-P(u_2,\ldots,u_r)\Bigr).\eqno(2.1.2)$$
This formula defines a partner for any mould $P\in ARI$, but in the case 
of polynomial-valued  moulds it corresponds to (2.1.1) in the sense that if 
$P=ma(p)$, then $P'=ma(p')$. 

Recall that the $push$-operator on a mould is an operator of order $r+1$ in
depth $r$ defined by
$$push(P)(u_1,\ldots,u_r)= P(-u_1-\cdots-u_r,u_1,\ldots,u_{r-1}),$$
and that a mould $P$ is said to be {\it push-invariant} if $P=push(P)$.
We say that a polynomial $p\in {\rm Lie}[a,b]$ is push-invariant
if $ma(p)$ is.
\vskip .2cm
\noindent {\bf Lemma 2.1.1.} {\it Let $p\in {\rm Lie}[a,b]$ and 
$p'\in \Q\langle a,b\rangle$ be polynomials with no linear term in $a$,
and let $D$ denote the derivation of $\Q\langle a,b\rangle$ given by 
$a\mapsto p$, $b\mapsto p'$. Then $D([a,b])=0$ if and only if $p$ is 
push-invariant and $p'$ is its partner, in which case $p'\in {\rm Lie}[a,b]$,
i.e.~$D\in {\rm Der}^0{\rm Lie}[a,b]$.}
\vskip .2cm
\noindent Proof. It follows immediately from the formula (2.1.1) for $p'$ 
that $\partial_a(p')=0$, so we can set 
$P=ma(p)=ma\bigl(D(a)\bigr)$ and $P'=ma(p')=ma\bigl(D(b)
\bigr)$.  Using the fact that $ma$ is a Lie algebra morphism (see Appendix)
and the first identity of (1.3.2) we find that 
$$ma\bigl(D([a,b]\bigr)=ma\bigl([D(a),b]+[a,D(b)]\bigr)=[P,B]-dur(P'),\eqno(2.1.3)$$ 
where $B=ma(b)$ is the mould concentrated in depth 1 given by $B(u_1)=1$.
Note that the mould $[P,B]-dur(P')$ is equal to zero in depths $r\le 1$. 

Let us first assume that $p$ is push-invariant, or equivalently, $P$ is 
push-invariant. Then $P'$ is the partner of $P$ as given in 
(2.1.2).  We have
$$[P,B](u_1,\ldots,u_r)=P(u_1,\ldots,u_{r-1})-P(u_2,\ldots,u_r)\eqno(2.1.4)$$
and
$$dur(P')=P(u_2,\ldots,u_r)-P(u_2,\ldots,u_{r-1},-u_1-\cdots-u_{r-1}).\eqno(2.1.5)$$
Thus $[P,B]-dur(P')$ is given in depth $r>1$ by
$$P(u_1,\ldots,u_{r-1})-P(u_2,\ldots,u_{r-1},-u_1-\cdots-u_{r-1})=
\bigl(P-push^{-1}(P)\bigr)(u_1,\ldots,u_r),$$
but since $P$ is push-invariant, this is equal to zero, so by (2.1.3)
$D([a,b])=0$. 

Assume now that $D([a,b])=0$, i.e., $[P,B]=dur(P')$, i.e.,
$$P(u_1,\ldots,u_{r-1})-P(u_2,\ldots,u_r)=(u_1+\cdots+u_r)P'(u_1,\ldots,u_r).\eqno(2.1.6)$$
This actually functions as a defining equation for $P'$. But knowing that
$P'=ma(p')$ is a polynomial-valued mould, (2.1.6) implies that
$P(u_1,\ldots,u_{r-1})-P(u_2,\ldots,u_r)$ must vanish along the pole
$u_1+\cdots+u_r=0$, in other words when $u_r=-u_1-\cdots-u_{r-1}$, so we have
$$P(u_1,\ldots,u_{r-1})=P(u_2,\ldots,u_{r-1},-u_1-\cdots-u_{r-1}).
\eqno(2.1.7)$$
As noted above, the right-hand side of (2.1.7) is nothing other than 
$push^{-1}(P)$, so (2.1.7) shows that $P$ is push-invariant. Furthermore, we can
substitute (2.1.7) into the left-hand side of (2.1.6) to find a new
defining equation for $P'$:
$$P'(u_1,\ldots,u_r)={{1}\over{u_1+\cdots+u_r}}\Bigl(P(u_2,\ldots,u_{r-1},
-u_1-\cdots-u_{r-1})-P(u_2,\ldots,u_r)\Bigr),\eqno(2.1.8)$$
but this coincides with (2.1.2), showing that $P'$ is the partner of 
$P$.  To conclude the proof, it remains only to show that in the case where
$p\in {\rm Lie}[a,b]$ is push-invariant, then its partner $p'$ is also a 
Lie polynomial, which is shown in [S2], Theorem 2.1 (i).
\hfill{$\diamondsuit$} 
\vskip .4cm
\noindent {\bf Lemma 2.1.2} {\it Let $\tilde D$ be a derivation of the Lie
subalgebra ${\rm Lie}[t_{01},t_{02}]\subset {\rm Lie}[a,b]$. Then
\vskip .1cm
\noindent (i) there exists a unique derivation $D$ of $\Q\langle a,b\rangle$
having the following two properties:
\vskip .1cm
(i.1) $D(t_{02})=\tilde D(t_{02})$;
\vskip .1cm
(i.2) $D(b)$ is the partner of $D(a)$.
\vskip .1cm
\noindent (ii) If $\tilde D(t_{12})=0$ and $D(a)$ is push-invariant, then 
$D$ is the unique extension of $\tilde D$ to all of ${\rm Lie}[a,b]$, so
in particular $D\in {\rm Der}^0{\rm Lie}[a,b]$.}
%
\vskip .2cm
\noindent Proof. (i) Let $T=\tilde D(t_{02})$, and write $T=\sum_{n\ge w} 
T_n$ for its homogeneous parts of weight $n$, where the weight is the degree 
as a polynomial in $a$ and $b$, and $w$ is the minimal weight occurring in 
$T$. We will construct a derivation $D$ satisfying $D(t_{02})=\tilde D(t_{02})$
by solving for $D(a)$ weight by weight via the equality
$$\eqalign{T&=D\bigl(Ber_{-b}(a)\bigr)\cr
&=D\bigl(a+{{1}\over{2}}[b,a]+
{{1}\over{12}}[b,[b,a]]-{{1}\over{720}}[b,[b,[b,[b,a]]]]+\cdots\bigr)\cr
&=D(a)+{{1}\over{2}}[D(b),a]+{{1}\over{2}}[b,D(a)]+{{1}\over{12}}[D(b),[b,a]]-{{1}\over{720}}[D(b),[b,[b,[b,a]]]
-{{1}\over{720}}[b,[D(b),[b,[b,a]]]]\cr
&\qquad\qquad\qquad\qquad-{{1}\over{720}}[b,[b,[D(b),[b,a]]]]+
\cdots.}\eqno(2.1.9)$$
We construct $D(a)$ by solving (2.1.9) in successive weights starting with $w$.
We start by setting $D(a)_w=T_w$ since $D(a)$ is the only term in (2.1.9) 
which can contribute to the lowest weight part $T_w$.
Let $D(b)_w$ be the partner of $D(a)_w$.  We then continue to solve the 
successive weight parts of (2.1.9) for $D(a)$ in terms of $T$ and the
previously determined lower weight 
parts of $D(a)$ and $D(b)$.  For instance 
the next few steps after weight $w$ are given by
$$\eqalign{D(a)_{w+1}&=T_{w+1}-{{1}\over{2}}[D(b)_w,a]-{{1}\over{2}}[b,D(a)_w],\cr
D(a)_{w+2}&=T_{w+2}-{{1}\over{2}}[D(b)_{w+1},a]-{{1}\over{2}}[b,D(a)_{w+1}]-{{1}\over{12}}[D(b)_w,[b,a]]\cr
&\qquad\qquad -{{1}\over{12}}[b,[D(b)_w,a]]-{{1}\over{2}}[b,[b,D(a)_w]],\cr
D(a)_{w+3}&=T_{w+3}-{{1}\over{2}}[D(b)_{w+2},a]-{{1}\over{2}}[b,D(a)_{w+2}]
 -{{1}\over{12}}[D(b)_{w+1},[b,a]]\cr
&\qquad\qquad -{{1}\over{12}}[b,[D(b)_{w+1},a]]
 -{{1}\over{12}}[b,[b,D(a)_{w+1}]]\cr
D(a)_{w+4}&=T_{w+4}-{{1}\over{2}}[D(b)_{w+3},a]-{{1}\over{2}}[b,D(a)_{w+3}]
-{{1}\over{12}}[D(b)_{w+2},[b,a]]\cr
&\qquad\qquad -{{1}\over{12}}[b,[D(b)_{w+2},a]]
-{{1}\over{12}}[b,[b,D(a)_{w+2}]]
+{{1}\over{720}} [D(b)_w,[b,[b,[b,a]]]\cr
&\qquad\qquad+{{1}\over{720}}[b,[D(b)_w,[b,[b,a]]]] +{{1}\over{720}}[b,[b,[D(b)_w,[b,a]]]]\cr
&\qquad\qquad+{{1}\over{720}}[b,[b,[b,[D(b)_w,a]]]] +{{1}\over{720}}[b,[b,[b,[b,D(a)_w]]]].}$$
In this way we construct the unique Lie series $D(a)$ and its partner
$D(b)$ such that the derivation $D$ satisfies 
$D\bigl(Ber_{-b}(a)\bigr)=D(t_{02})=T=\tilde D(t_{02})$. We note that
$D(b)$ is not necessarily a Lie polynomial, and furthemore $D$ is not 
necessarily an extension of $\tilde D$ to all of ${\rm Lie}[a,b]$, 
because $D$ and $\tilde D$ may not agree on $t_{12}$. As an example, consider
the case where $\tilde D(t_{02})=t_{12}$ and $\tilde D(t_{12})=0$. The
above process to construct $D(a)$ shows that in minimal weight $w$ we have
$D(a)_2=t_{12}=[a,b]$. But by the formula (2.1.1) for the partner, we find that
$D(b)_2=-b^2$. The derivation of $\Q\langle a,b\rangle $ mapping 
$a\mapsto [a,b]$ and $b\mapsto -b^2$ does not annihilate $t_{12}$, so
$D$ is not an extension of $\tilde D$ to all of $\Q\langle a,b\rangle$.

To prove (ii), we first assume that $D(a)$ is push-invariant. Then 
since $D(b)$ is the partner of $D(a)$ by definition, Lemma 2.1.1 shows that
$D(b)$ is a Lie series and $D$ annihilates $t_{12}=[a,b]$. Thus $D$
is an extension of $\tilde D$ to a derivation of all of ${\rm Lie}[a,b]$.
For the uniqueness, suppose that $E$ is another derivation of ${\rm Lie}[a,b]$ 
that coincides with $\tilde D$ on $t_{02}$ and $t_{12}$.  The fact that
$E(t_{12})=E([a,b])=0$ shows that $E(a)$ and $E(b)$ are partners by Lemma 2.1.1.
But then $E$ satisfies (i.1) and (i.2), so it coincides with $D$.\hfill{$\diamondsuit$}
\vskip .5cm
\noindent {\bf 2.2. Proof of the main theorem.}
\vskip .2cm
For each $\psi\in \grt$, let $f(x,y)=\psi(x,-y)$.  Let
$A=Ad_{ari}(invpal)\cdot ma(f)$ as before, and $M=\Delta(A)$.
By Corollary 1.3.4, there exists a polynomial $m\in {\rm Lie}[C]$ 
such that
$$ma(m)=M=\Delta\Bigl(Ad_{ari}(invpal)\cdot ma(f)\Bigr).$$
Since by the same corollary $m$ is push-invariant, we see that by
Lemma 2.1.1 there exists a unique derivation $E_\psi\in \Der\,{\rm Lie}[a,b]$
such that $E_\psi(a)=m$, $E_\psi([a,b])=0$ and $E_\psi(b)\in {\rm Lie}[C]$,
namely the one such that $E_\psi(b)$ is the partner
of $E_\psi(a)$.  The main result we need about this derivation is the following.
\vskip .2cm
\noindent {\bf Proposition 2.2.1.} {\it The derivation $E_\psi$ satisfies
$$E_\psi(t_{02})=[\psi(t_{02},t_{12}),t_{02}].\eqno(2.2.1)$$}
Using this, we can easily prove the main theorem.  Since 
$t_{12}=[a,b]$, we have $E_\psi(t_{12})=0$, so Proposition 2.2.1 shows that
$E_\psi$ restricts to a derivation $\tilde E_\psi$ on the Lie subalgebra 
${\rm Lie}[t_{02},t_{12}]$, where it coincides with the restriction
$\tilde D_\psi$ of Enriquez' derivation $D_\psi$ given in 
(1.2.2).  Furthermore, since $E_\psi(t_{12})=0$ and $E_\psi(a)=m$ 
is push-invariant, we are in the situation of Lemma 2.1.2 (ii), so $E_\psi$ is 
the unique extension of $\tilde E_\psi$ to all of 
${\rm Lie}[a,b]$.  But Enriquez' derivation $D_\psi$ is an extension of
$\tilde D_\psi$ to all of ${\rm Lie}[a,b]$, and it also satisfies 
$D_\psi(t_{12})=0$, so by Lemma 2.1.1, $D_\psi(a)=\alpha_+=\gamma_+(\psi)$ 
is push-invariant; thus by Lemma 2.1.2 (ii) $D_\psi$ is the 
unique extension of $\tilde D_\psi$ to all of ${\rm Lie}[a,b]$. 
Thus, since $\tilde E_\psi= \tilde D_\psi$, we must have $E_\psi=D_\psi$, and 
in particular $E_\psi(a)=m=D_\psi(a)=\gamma_+(\psi)$.  Taking $ma$ of both sides 
yields the desired equality (1.3.5).\hfill{$\diamondsuit$}
\vskip .5cm
\noindent {\bf 3. Proof of Proposition 2.2.1}
\vskip .3cm
\noindent {\bf 3.1. Mould theoretic derivations}
\vskip .3cm
We begin by defining a mould-theoretic derivation ${\cal E}_\psi$ on 
$ARI_{lu}$ for each $\psi\in\grt$ as follows.
\vskip .3cm
\noindent {\bf Definition.} For any mould $P$, let $Darit(P)$ be the operator 
on moulds defined by
$$Darit(P)=-dar\Bigl(arit\bigl(\Delta^{-1}(P)\bigr)-ad\bigl(\Delta^{-1}(P)\bigr)\Bigr)\circ dar^{-1}.\eqno(3.1.1)$$
Then for all $P$, $Darit(P)$ is a derivation of $ARI_{lu}$,
since $arit(P)$ and $ad(P)$ are both derivations and $dar$ is an automorphism.

Let $\psi\in\grt$. We use the notation of (1.3.8), and set
$${\cal E}_\psi=Darit(M).\eqno(3.1.2)$$

Recall that $ARI$ denotes the vector space of rational-valued moulds with 
constant term $0$.  
Let $ARI^a$ denote the vector space obtained by adding a single generator
$a$ to the vector space $ARI$, and let $ARI^a_{lu}$ be the Lie algebra
formed by extending the $lu$-bracket to $ARI^a$ via the relation
$$[Q,a]=dur(Q)\eqno(3.1.3)$$
for every $Q\in ARI_{lu}$.  Recall from (1.3.2) that this equality holds in
the polynomial sense if $Q$ is a polynomial-valued mould; in
other words, (1.3.3) extends to an injective Lie algebra morphism
$ma:{\rm Lie}[a,b]\rightarrow ARI^a_{lu}$ by formally setting $ma(a)=a$.

The Lie algebra $ARI_{lu}$ forms a Lie ideal of $ARI^a_{lu}$,
i.e., there is an exact sequence of Lie algebras
$$0\rightarrow ARI_{lu}\rightarrow ARI^a_{lu}\rightarrow \Q a \rightarrow 0.$$
We say that a derivation (resp. automorphism) of $ARI_{lu}$ {\it extends to
$a$} if there is a derivation (resp. automorphism)
of $ARI^a_{lu}$ that restricts to the given one on the Lie subalgebra
$ARI_{lu}$.  To check whether a given derivation (resp. automorphism)
extends to $a$, it suffices to check that relation (3.1.3) is respected.

Recall that $B=ma(b)$ is the mould concentrated in depth $1$ given by
$B(u_1)=1$.  Let us write $B_i$, $i\ge 0$, for the mould concentrated
in depth $1$ given by $B_i(u_1)=u_1^i$. In particular $B_0=B=ma(b)$,
and $B_1(u_1)=u_1$, so $B_1=ma([b,a])$.
\vskip .3cm
\noindent {\bf Lemma 3.1.1.} {\it (i) The automorphism $dar$ extends to $a$
taking the value $dar(a)=a$;
\vskip .1cm
\noindent (ii) The derivation $dur$ extends to $a$ taking the value $dur(a)=0$;
\vskip .1cm
\noindent (iii) For all $P\in ARI$, the derivation
$arit(P)$ of $ARI_{lu}$ extends to $a$, taking the value $arit(P)\cdot a=0$.
\vskip .1cm
\noindent (iv) For all $P\in ARI$, the derivation $Darit(P)$
of $ARI_{lu}$ extends to $a$, with $Darit(P)\cdot a=P$. Furthermore,
$Darit(P)\cdot B_1=0$.}
\vskip .2cm
\noindent {\bf Proof.} Since $dar$ is an automorphism, to check (3.1.3) we 
write
$$[dar(Q),dar(a)]=[dar(Q),a]=dur\bigl(dar(Q)\bigr).$$ 
But it is obvious from their definitions that $dur$ and $dar$ commute, 
so this is indeed equal to $dar\bigl(dur(Q)\bigr)$. This proves (i).
We check (3.1.3) for (ii) similarly. Because $dur(a)=0$ and
$dur$ is a derivation, we have
$$dur([Q,a])=[dur(Q),a]=dur\bigl(dur(Q)\bigr).$$

For (iii), we have
$$arit(P)\cdot [Q,a]=[arit(P)\cdot Q,a]=dur\bigl(arit(P)\cdot Q)\bigr).$$
But as pointed out by Ecalle [E2] (cf. [S, Lemma 4.2.2] for details), 
$arit(P)$ commutes with $dur$ for all $P$, which proves the result. 

For (iv), the calculation to check that (3.1.3) is respected is a little
more complicated.  Let $Q\in ARI$. Again using the
commutation of $arit(P)$ with $dur$, as well as that of $dar$ and $dur$,
we compute
$$\eqalign{Darit&(P)\cdot [Q,a]=\bigl[Darit(P)(Q),a\bigr]+
\bigl[Q,Darit(P)(a)\bigr]\cr
&=dur\bigl(Darit(P)\cdot Q\bigr)+[Q,P]\cr
&=-dur\Biggl(dar\Bigl(arit\bigl(\Delta^{-1}(P)\bigr)\cdot dar^{-1}(Q)-\bigl[\Delta^{-1}(P),dar^{-1}(Q)\bigr)\Bigr]\Biggr)+[Q,P]\cr
&=-dur\Biggl(dar\Bigl(arit\bigl(\Delta^{-1}(P)\bigr)\cdot dar^{-1}(Q)\Bigr)\Biggr)-dur\Bigl(\bigl[Q,dur^{-1}(P)\bigr]\Bigr)+[Q,P]\cr
&=-dar\Biggl(dur\Bigl(arit\bigl(\Delta^{-1}(P)\bigr)\cdot dar^{-1}(Q)\Bigr)\Biggr)-\bigl[[Q,N],a\bigr]+\bigl[Q,[N,a]\bigr]\cr
&\qquad\qquad\qquad\qquad\qquad\qquad{\rm with\ } N=dur^{-1}P,\ {\rm i.e.,}\ P=[N,a]\cr}$$
$$\eqalign{\qquad\qquad 
&=-dar\Bigl(arit\bigl(\Delta^{-1}(P)\bigr)\cdot dur\,dar^{-1}(Q)\Bigr)
-\bigl[[Q,a],N\bigr]\ \ {\rm by\ Jacobi}\cr 
&=-dar\Bigl(arit\bigl(\Delta^{-1}(P)\bigr)\cdot dar^{-1}\,dur(Q)\Bigr)
-\bigl[dur(Q),dur^{-1}P\bigr]\cr
&=-dar\Bigl(arit\bigl(\Delta^{-1}(P)\bigr)\cdot dar^{-1}\,dur(Q)\Bigr)
-dar\Bigl(\bigl[dar^{-1}dur(Q),dar^{-1}dur^{-1}(P)\bigr]\Bigr)\cr
&=-dar\Bigl(arit\bigl(\Delta^{-1}(P)\bigr)\cdot dar^{-1}\,dur(Q)\Bigr)
+dar\Bigl(\bigl[\Delta^{-1}(P),dar^{-1}dur(Q)\bigr]\Bigr)\cr
&=Darit(P)\cdot dur(Q).}$$
\vskip -.2cm
This proves the first statement of (iv).  For the second statement, we note that
$dar^{-1}(B_1)=B$.  Set $R=\Delta^{-1}(P)$.  We compute
$$\eqalign{Darit(P)\cdot B_1&=
-dar\bigl(arit(R)\cdot B\bigr)+dar\bigl([R,B]\bigr)\cr
&=-u_1\cdots u_r \bigl(R(u_1,\ldots,u_{r-1})-R(u_2,\ldots,u_r)\bigr)\cr
&\qquad -u_1\cdots u_r\bigl(R(u_2,\ldots,u_r)-R(u_1,\ldots,u_{r-1})\bigr)\cr
&=0.}$$
This concludes the proof of Lemma 3.1.1.\hfill{$\diamondsuit$}
\vskip .3cm
We consider by default that $a$ is alternal and polynomial.
Let $(ARI^a_{lu})^{pol,al}$ denote the Lie subalgebra of alternal polynomial
moulds of $ARI^a_{lu}$. Then $ARI^{pol,al}_{lu}$ is a Lie ideal of
$ARI^a_{lu}$ and we have the Lie algebra isomorphism 
$$L[C]{\Bbb o}\Q a\simeq {\rm Lie}[a,b]\buildrel{ma}\over\longrightarrow
\bigl(ARI^a_{lu}\bigr)^{pol,al} \simeq ARI^{pol,al}_{lu}
{\Bbb o}\Q a.\eqno(3.1.4)$$
\vskip .3cm
\noindent {\bf Lemma 3.1.2.} {\it Suppose that $P\in ARI$ is a mould such that 
$Darit(P)$ preserves the Lie subalgebra $(ARI^a_{lu})^{pol,al}$ of
$ARI^a_{lu}$.
Then there exists a derivation $E_P\in \Der\,{\rm Lie}[a,b]$ that corresponds to $Darit(P)$ restricted to $(ARI^a_{lu})^{pol,al}$, in the sense that
$$ma\bigl(E_P(f)\bigr)=Darit(P)\bigl(ma(f)\bigr)\ \ \ {\rm for\ all}\ 
f\in {\rm Lie}[a,b].$$
The derivation $E_P$ has the property that the values $E_P(a)$ and
$E_P(b)$ lie in ${\rm Lie}[C]$.}
\vskip .2cm
\noindent Proof. By the isomorphism (3.1.4), every mould
$P\in (ARI^a_{lu})^{pol,al}$ has a unique preimage in ${\rm Lie}[a,b]$
under $ma$: we write $p=ma^{-1}(P)$.  Recall that $B=ma(b)$. By assumption, 
$P$ is an alternal polynomial-valued mould, and so is $Darit(P)\cdot B$ 
since $P$ preserves such moulds. Thus we can define $E_P$ by setting
$E_P(a)=ma^{-1}(P)$, $E_P(b)=ma^{-1}\bigl(Darit(P)\cdot B\bigr)$. In particular
this means that the monomial $a$ does not appear in the polynomials
$E_P(a)$ and $E_P(b)$.  \hfill{$\diamondsuit$}
\vskip .3cm
\noindent {\bf Lemma 3.1.3.} {\it Let $P$ be an alternal polynomial-valued 
mould.  Then $Darit(P)$ preserves $(ARI^a_{lu})^{pol,al}$ if and only if $P$ 
is push-invariant.}
\vskip .2cm
\noindent Proof. By the isomorphism (3.1.4), $(ARI^a_{lu})^{pol,al}$ is
generated as a Lie algebra under the $lu$ bracket by $ma(a)=a$ and $ma(b)=B$.
Since $Darit(P)\cdot a=P$ is alternal and polynomial-valued
by assumption, it suffices to determine when $Darit(P)\cdot B$
is alternal and polynomial.  Let $N=\Delta^{-1}P$,
and set $B_{-1}=dar^{-1}(B)$, so $B_{-1}$ is concentrated in depth 1 with
$B_{-1}(u_1)=1/u_1$. We compute
$$\eqalign{\bigl(Darit(P)\cdot &B\bigr)(u_1,\ldots,u_r)=
-dar\bigl(arit(N)\cdot B_{-1}-[N,B_{-1}]\bigr)(u_1,\ldots,u_r)\cr
&=-dar\bigl(arit(N)\cdot B_{-1}\bigr)(u_1,\ldots,u_r)
-dar\bigl([B_{-1},N]\bigr)(u_1,\ldots,u_r)\cr
&=-dar\Bigl(B_{-1}(u_1+\cdots+u_r)\bigl(N(u_1,\ldots,u_{r-1})-N(u_2,\ldots,u_r)\bigr)\Bigr)\cr
&\qquad -u_1\ldots u_r\bigl(B_{-1}(u_1)N(u_2,\ldots,u_r)+N(u_1,\ldots,u_{r-1})B_{-1}(u_r)\bigr)\cr
&=-u_1\cdots u_r(u_1+\cdots+u_r)^{-1}\bigl(N(u_1,\ldots,u_{r-1})-N(u_2,\ldots,u_r)\bigr)\cr
&\qquad -u_2\cdots u_rN(u_2,\ldots,u_r)+u_1\cdots u_{r-1}N(u_1,\ldots,u_{r-1})\cr
&={{1}\over{u_1+\cdots+u_r}}\bigl(P(u_1,\ldots,u_{r-1})-P(u_2,\ldots,u_r)\bigr).}$$
In order for this mould to be polynomial-valued, it is necessary and sufficient that the numerator should be zero when $u_r=-u_1-\cdots-u_{r-1}$, i.e., that
$$P(u_1,\ldots,u_{r-1})=P(u_2,\ldots,u_{r-1},-u_1-\cdots-u_{r-1}).\eqno(3.1.5)$$
But the right-hand term is equal to $push^{-1}(P)$, so this condition is
equivalent to the push-invariance of $P$.\hfill{$\diamondsuit$}
\vskip .3cm
\noindent {\bf Corollary 3.1.4.} {\it The derivation $E_\psi$ defined in 
section 2.2 is equal to the derivation $E_M$ associated to $Darit(M)$ as in 
Lemma 3.1.2.}
\vskip .2cm
\noindent Proof. Since $M$ is push-invariant by Corollary 1.3.4, 
$Darit(M)$ preserves $(ARI^a_{lu})^{pol,al}$ by Lemma 3.1.3. Thus
we are in the situation of Lemma 3.1.2, so there exists
a derivation $E_M$ of ${\rm Lie}[a,b]$ such that $E_M(a)=m$ with
$ma(m)=M$. Furthermore, setting $B_1=ma([b,a])$, we know that
$Darit(M)\cdot B_1=0$ by Lemma 3.1.1 (iv), and therefore by Lemma 3.1.2, we
have $E_M([b,a])=E_M([a,b])=0$. Thus the derivation $E_M$ of ${\rm Lie}[a,b]$
agrees with $E_\psi$ on $a$ and on $[a,b]$, so since furthermore
$E_M(b)\in {\rm Lie}[a,b]\ominus {\rm Lie}[a]$, they are equal.
\hfill{$\diamondsuit$}
\vskip .3cm
This result means that we can now use mould theoretic methods to study
$Darit(M)$ in order to prove Proposition 2.2.1.
\vskip .4cm
\noindent {\bf 3.2. The $\Delta$-operator}
\vskip .3cm
Let us define a new Lie bracket, the $Dari$-bracket, on $ARI$ by 
$$Dari(P,Q)=Darit(P)\cdot Q-Darit(Q)\cdot P,$$
where $Darit(P)$ is the $lu$-derivation defined in (3.1.1).  Let
$ARI_{Dari}$ denote the Lie algebra obtained by equipping $ARI$ with
this Lie bracket.
\vskip .3cm
\noindent {\bf Proposition 3.2.1.} {\it The operator $\Delta$ is a Lie
algebra isomorphism from $ARI_{ari}$ to $ARI_{Dari}$. } 
\vskip .2cm
\noindent Proof. Certainly
$\Delta$ is a vector space isomorphism from $ARI_{ari}$ to $ARI_{Dari}$
since it is an invertible operator on moulds. To prove that it is a
Lie algebra isomorphism, we need to show the Lie bracket identity
$\Delta\bigl(ari(P,Q)\bigr)=Dari\bigl(\Delta P,\Delta Q\bigr)$, or
equivalently,
$$Dari(P,Q)=\Delta\bigl(ari(\Delta^{-1}P,\Delta^{-1}Q)\bigr)\eqno(3.2.1)$$
for all moulds $P,Q\in ARI$. But indeed, we have
$$\eqalign{Dari(P,Q)&=Darit(P)\cdot Q-Darit(Q)\cdot P\cr
&=-\bigl(dar\circ arit(\Delta^{-1}P)\circ dar^{-1}\bigr)\cdot Q+
\bigl(dar\circ ad(\Delta^{-1}P)\circ dar^{-1}\bigr)\cdot Q\cr
&\qquad\qquad +\bigl(dar\circ arit(\Delta^{-1}Q)\circ dar^{-1}\bigr)\cdot P
-\bigl(dar\circ ad(\Delta^{-1}Q)\circ dar^{-1}\bigr)\cdot P\cr
&=-\bigl(\Delta\circ arit(\Delta^{-1}P)\circ \Delta^{-1}\bigr)\cdot Q
+\bigl(\Delta\circ arit(\Delta^{-1}Q)\circ \Delta^{-1}\bigr)\cdot P\cr
&\qquad\qquad +\bigl(dar\circ ad(\Delta^{-1}P)\circ dar^{-1}\bigr)\cdot Q
-\bigl(dar\circ ad(\Delta^{-1}Q)\circ dar^{-1}\bigr)\cdot P\cr
&=-\bigl(\Delta\circ arit(\Delta^{-1}P)\circ \Delta^{-1}\bigr)\cdot Q
+\bigl(\Delta\circ arit(\Delta^{-1}Q)\circ \Delta^{-1}\bigr)\cdot P\cr
&\qquad\qquad +dar\bigl([\Delta^{-1}(P),dar^{-1}Q]\bigr)-
dar\bigl([\Delta^{-1}(P),dar^{-1}P]\bigr)\cr
&=\Delta\Bigl(-arit(\Delta^{-1}P\cdot \Delta^{-1}Q+arit(\Delta^{-1}Q)\cdot
\Delta^{-1}P\cr
&\qquad\qquad +dur^{-1}\bigl([\Delta^{-1}P,dar^{-1}Q]+[dar^{-1}P,\Delta^{-1}Q]\bigr)\Bigr)\cr
&=\Delta\Bigl(-arit(\Delta^{-1}P\cdot \Delta^{-1}Q+arit(\Delta^{-1}Q)\cdot
\Delta^{-1}P\cr
&\qquad\qquad +dur^{-1}\bigl([\Delta^{-1}P,dur\Delta^{-1}Q]+[dur\Delta^{-1}P,\Delta^{-1}Q]\bigr)\Bigr)\cr
}$$
$$\eqalign{
&=\Delta\Bigl(-arit(\Delta^{-1}P\cdot \Delta^{-1}Q+arit(\Delta^{-1}Q)\cdot
\Delta^{-1}P\cr
&\qquad\qquad +dur^{-1}dur\bigl([\Delta^{-1}P,\Delta^{-1}Q]\bigr)\Bigr)\cr
&=\Delta\Bigl(-arit(\Delta^{-1}P\cdot \Delta^{-1}Q+arit(\Delta^{-1}Q)\cdot
\Delta^{-1}P+[\Delta^{-1}P,\Delta^{-1}Q]\Bigr)\cr
&=\Delta\bigl(ari(\Delta^{-1}P,\Delta^{-1}Q)\bigr),}$$
which proves the desired identity.  \hfill{$\diamondsuit$}
\vfill\eject
Let us now define the group $GARI_{Dgari}$. We start by defining
the exponential map $exp_{Dari}:ARI_{Dari}\rightarrow GARI$ by
$$exp_{Dari}(P)=1+\sum_{n\ge 1} {{1}\over{n!}}Darit(P)^{n-1}(P),\eqno(3.2.2)$$
which for all $P\in ARI$ satisfies the equality
$$exp\bigl(Darit(P)\bigr)(a)=exp_{Dari}(P).\eqno(3.2.3)$$
This map is easily seen to be invertible, since for any $Q\in GARI$ we can
recover $P$ such that $exp_{Dari}(P)=Q$ recursively depth by depth. Let
$log_{Dari}$ denote the inverse of $exp_{Dari}$.
For each $P\in GARI$, we then define an automorphism 
$Dgarit(P)\in Aut\,ARI_{lu}$ by 
$$Dgarit(P)=Dgarit\Bigl(exp_{Dari}\bigl(log_{Dari}(P)\bigr)\Bigr)=
exp\Bigl(Darit\bigl(log_{Dari}(P)\bigr)\Bigr).$$
Finally, we define the 
multiplication $Dgari$ on $GARI$ by
$$\eqalign{Dgari(P,Q)&=exp_{Dari}\bigl(ch_{Dari}(log_{Dari}(P),log_{Dari}(Q))\bigr)\cr
&=exp\bigl(Darit(log_{Dari}(P))\bigr)\circ exp\bigl(Darit(log_{Dari}(Q))\bigr)\cdot a\cr
&=Dgarit(P)\circ Dgarit(Q)\cdot a\cr
&=Dgarit(P)\cdot Q,}$$
where $ch_{Dari}$ denotes the Campbell-Hausdorff law on $ARI_{Dari}$.
We obtain the following commutative diagram, analogous to Ecalle's diagram (A.18) 
(cf. Appendix):
$$\xymatrix{ARI_{Dari}\ar[r]^{exp_{Dari}\ }\ar[d]_{Darit}&GARI_{Dgari}
\ar[d]^{Dgarit}\\
{\rm \Der}\,ARI_{lu}\ar[r]^{exp}&{\rm Aut}\,ARI_{lu}.}\eqno(3.2.4)$$
\vskip .3cm
\noindent {\bf Lemma 3.2.2.} {\it For any mould $P\in GARI$, the
automorphism $Dgarit(P)$ of $ARI_{lu}$ extends to an automorphism of the 
Lie algebra $ARI^a_{lu}$ with the following properties:
\vskip .1cm
i) its value on $a$ is given by
$$Dgarit(P)\cdot a=a-1+P\in ARI^a;\eqno(3.2.5)$$
\vskip .1cm
ii) we have $Dgarit(P)\cdot B_1=B_1$.}
\vskip .2cm
\noindent Proof. Let $Q=log_{Dari}(P)\in ARI$. We saw in Lemma 3.1.1 (iv)
that $Darit(Q)$ extends to $ARI^a_{lu}$ with $Darit(Q)\cdot a=Q$. 
By diagram (3.2.4), we have
$$\eqalign{Dgarit(P)\cdot a&=Dgarit\bigl(exp_{Dari}(Q)\bigr)\cdot a\cr
&=exp\bigl(Darit(Q)\bigr)\cdot a\cr
&=a+Darit(Q)\cdot a + {{1}\over{2}}Darit(Q)^2\cdot a+\cdots\cr
&=a+Q + {{1}\over{2}}Darit(Q)\cdot Q+\cdots\cr
&=a-1+exp_{Dari}(Q)\ \ {\rm by\ (3.2.2)}\cr
&=a-1+P.}$$
The second statement follows immediately from the fact that
$Darit(Q)\cdot B_1=0$ for all $Q\in ARI$ shown in Lemma 3.1.1 (iv).\hfill{$\diamondsuit$}
\vskip .3cm
Finally, we set $\Delta^*=exp_{Dari}\circ \Delta\circ log_{ari}$, to obtain
the commutative diagram of isomorphisms
$$ARI_{ari}\buildrel{\Delta}\over\longrightarrow\ \ \ ARI_{Dari}$$
$$exp_{ari}\,\downarrow\ \ \ \qquad\qquad\downarrow\,exp_{Dari}$$
$$GARI_{gari}\buildrel{\Delta^*}\over\rightarrow GARI_{Dgari},\eqno(3.2.6)$$
which will play a special role in the proof of Proposition 2.2.1. Indeed,  
the key result in our proof Proposition 2.2.1 is an explicit formula for 
the map $\Delta^*$. In order to formulate it, we first define the 
{\it $mu$-dilator} of a mould, introduced by Ecalle in [E2].
\vskip .2cm
\noindent {\bf Definition.} Let $P\in GARI$.  Then the $mu$-dilator of
$P$, denoted $duP$, is defined by 
$$duP=P^{-1}\,dur(P).\eqno(3.2.7)$$
Ecalle writes this in the equivalent form $dur(P)=P\,duP$,
and by (3.1.3),
this means that $[P,a]=Pa-aP=P\,duP=P$, whch multiplying by $P^{-1}$, gives
us the useful formulation\footnote{$^6$}{We are grateful to B. Enriquez for spotting this 
enlightening interpretation of the $mu$-dilator, which cannot even be
stated meaningfully for general moulds unless $a$ is added to $ARI$.}
$$P^{-1}aP=a-duP.\eqno(3.2.8)$$
\noindent {\bf Proposition 3.2.3.} {\it The isomorphism
$$\Delta^*:GARI_{gari}\rightarrow GARI_{Dgari}$$
in diagram (3.2.6) is explicitly given by the formula
$$\Delta^*(Q)=1-dar\bigl(du\,inv_{gari}(Q)\bigr).\eqno(3.2.9)$$\par}
\noindent Proof.  Let $Q\in GARI$, and set $P=log_{ari}(Q)$. 
Let $R=exp_{ari}(-P)$. By Lemma A.1 from the Appendix, the derivation
$-arit(P)+ad(P)$ extends to $a$ taking the value $[a,P]$ on $a$, and we have
$$exp\bigl(-arit(P)+ad(P)\bigr)\cdot a=R^{-1}\,a\,R.\eqno(3.2.10)$$
By (3.1.1), we have
$$exp\Bigl(Darit\bigl(\Delta(P)\bigr)\Bigr)=
dar\circ exp\bigl(-arit(P)+ad(P)\bigr)\circ dar^{-1}.$$
Recall that $dar(a)=a$ by Lemma 3.1.1 (i), and $dar$ is
an automorphism of $ARI^a_{lu}$; in particular $du$ commutes with $dar$.  
Thus we have
$$\eqalign{exp\Bigl(Darit\bigl(\Delta(P)\bigr)\Bigr)\cdot a&=dar\circ 
exp\bigl(-arit(P)+ad(P)\bigr)\cdot a\cr
&=dar(R^{-1}\,a\,R)\ \ \ {\rm by\ Lemma\ A.1}\cr
&=dar(R)^{-1}\,a\,dar(R)\cr
&=a-du\bigl(dar(R)\bigr)\ \ {\rm by\ (3.2.8)}\cr
&=a-dar\bigl(duR\bigr).}\eqno(3.2.11)$$
Now, using $P=log_{ari}(Q)$, we compute 
$$\eqalign{\Delta^*(Q)&=1-a+Dgarit\bigl(\Delta^*(Q)\bigr)\cdot a
\ \ {\rm by\ (3.2.5)}\cr
&=1-a+Dgarit\Bigl(exp_{Dari}\bigl(\Delta(log_{ari}(Q))\bigr)\Bigr)\cdot a
\ \ {\rm by\ (3.2.6)}\cr
&=1-a+Dgarit\Bigl(exp_{Dari}\bigl(\Delta(P)\bigr)\Bigr)\cdot a\cr
&=1-a+exp\Bigl(Darit\bigl(\Delta(P)\bigr)\Bigr)\cdot a
\ \ {\rm by\ (3.2.4)}\cr
&=1-dar\bigl(du\,exp_{ari}(-P)\bigr)
\ {\rm by\ (3.2.11)}\cr
&=1-dar\bigl(du\,inv_{gari}(Q)\bigr).}\eqno(3.2.12)$$
This proves the proposition.  \hfill{$\diamondsuit$}
\vskip .3cm
\noindent {\bf Corollary.} {\it We have the
identity $$\Delta^*(invpal)=ma\bigl(1-a+Ber_{-b}(a)\bigr).\eqno(3.2.13)$$}
\noindent Proof. Applying (3.2.9) to $Q=invpal=inv_{gari}(pal)$, 
we find 
$$\Delta^*(invpal)=1-
dar\bigl(dupal\bigr),\eqno(3.2.14)$$
where $dupal$ is the $mu$-dilator of $pal$ given in (1.3.3),
discovered by Ecalle. Comparing the elementary mould identity
$$ma\Bigl(ad(-b)^r(-a)\Bigr)=\sum_{j=0}^{r-1} (-1)^j\Bigl({{r-1}\atop{j}}\Bigr)u_{j+1}$$
with (1.3.3) shows that $dar(dupal)$ is given in depth $r\ge 1$ by
$$dar(dupal)(u_1,\ldots,u_r)
={{B_r}\over{r!}}\sum_{j=0}^{r-1}(-1)^j\Bigl({{r-1}\atop{j}}\Bigr)u_{j+1}
={{B_r}\over{r!}}ma\Bigl(ad(-b)^r(-a)\Bigr).$$
Since the constant term of $dar\bigl(dupal\bigr)(\emptyset)$ is $0$, 
this yields 
$$dar\bigl(dupal\bigr)=ma\Bigl(Ber_{-b}(-a)+a\Bigr)=ma\bigl(a-Ber_{-b}(a)\bigr),$$
so (3.2.14) implies the desired identity (3.2.13).  \hfill{$\diamondsuit$}
\vskip .5cm
\noindent {\bf 3.3. Proof of Proposition 2.2.1}
\vskip .2cm
Let $\psi\in \grt$. We return to the notation of (1.3.8).
By Corollary 3.1.4, we have a derivation $E_M=E_\psi\in
\Der\,{\rm Lie}[a,b]$ obtained by restricting the derivation
${\cal E}_\psi=Darit(M)$ to the Lie subalgebra of $ARI^a_{lu}$ generated by
$a$ and $B=ma(b)$, which is precisely $(ARI^a_{lu})^{pol,al}$, and transporting the derivation to the isomorphic 
space ${\rm Lie}[a,b]$.  The purpose of this section is to prove (2.2.1), i.e.,
$$E_\psi(t_{02})=[\psi(t_{02},t_{12}),t_{02}].$$
The main point is the following result decomposing $Darit(M)$ into three 
factors; a derivation conjugated by an automorphism.  We note that although the
values of the derivation and the automorphism in Proposition 3.3.1 on $a$ are 
polynomial-valued moulds, this is false for their values on $B=ma(b)$, which 
means that this decomposition is a result which cannot be stated in the
power-series situation of ${\rm Lie}[a,b]$; the framework of mould
theory admitting denominators is crucial here.
\vskip .3cm
\noindent {\bf Proposition 3.3.1.} {\it We have the following identity of 
derivations:
$$Darit\Bigl(\Delta\bigl(Ad_{ari}(invpal)\cdot F\bigr)\Bigr) = $$
$$Dgarit\bigl(\Delta^*(invpal)\bigr)\circ
Darit\bigl(\Delta(F)\bigr)\circ Dgarit\bigl(\Delta^*(invpal)\bigr)^{-1}.\eqno(3.3.1)$$}

\noindent Proof. We use two standard facts about Lie algebras and their
exponentials. Firstly,
for any exponential morphism $exp:\frak{g}\rightarrow G$ mapping a Lie algebra 
to its associated group, the natural adjoint action of $G$ on $\frak{g}$, 
denoted $Ad_\frak{g}(exp(g))\cdot h$, satisfies
$$exp\Bigl(Ad_\frak{g}\bigl(exp(g)\bigr)\cdot h\Bigr)=
Ad_G\bigl(exp(g)\bigr)\bigl(exp(h)\bigr)=exp(g)*_G exp(h)*_G exp(g)^{-1},
\eqno(3.3.2)$$
where $*_G$ denotes the multiplication in $G$, defined by
$$exp(g)*_G exp(h)=exp\bigl(ch_\frak{g}(g,h)\bigr)\eqno(3.3.3)$$
where $ch_\frak{g}$ denotes the Campbell-Hausdorff law on $\frak{g}$.

Secondly, if $\Delta:\frak{g}\rightarrow\frak{h}$ is an isomorphism of Lie 
algebras, then the following diagram commutes:
$$\xymatrix{\frak{g}\ar[r]^{\Delta}\ar[d]_{Ad_{\frak{g}}
\bigl(exp_\frak{g}(g)\bigr)} &\frak{h}
\ar[d]^{Ad_{\frak{h}}\Bigl(exp_{\frak{h}}\bigl(\Delta(g)\bigr)\Bigr)}\\
\frak{g}\ar[r]^{\Delta}&\frak{h}.}\eqno(3.3.4)$$

To prove (3.3.1), we start by taking the exponential of both sides. Let
$lipal=log_{ari}(invpal)$. We start
with the left-hand side and compute
$$\eqalign{&exp\Biggl(Darit\Bigl(\Delta\bigl(Ad_{ari}(invpal)\cdot F\bigr)\Bigr)\Biggr) = 
exp\Biggl(Darit\Bigl(\Delta\bigl(Ad_{ari}(exp_{ari}(lipal))\cdot F\bigr)\Bigr)\Biggr)\cr
&\ =exp\Biggl(Darit\Bigl(Ad_{Dari}\bigl(exp_{Dari}(\Delta lipal)\bigr)\cdot \Delta(F)\bigr)\Bigr)\Biggr)\cr
&\ =Dgarit\Biggl(exp_{Dari}\Bigl(Ad_{Dari}\bigl(exp_{Dari}(\Delta lipal)\bigr)\cdot \Delta(F)\Bigr)\Biggr)\qquad\qquad\qquad\qquad\qquad\quad {\rm (3.3.5)}\cr
&\ =Dgarit\Bigl(exp_{Dari}\bigl(\Delta lipal\bigr)\Bigr)\circ Dgarit\Bigl(exp_{Dari}\bigl(\Delta(F)\bigr)\Bigr)\circ Dgarit\Bigl(exp_{Dari}\bigl(\Delta lipal\bigr)\Bigr)^{-1}\cr
&\ =Dgarit\bigl(\Delta^*(invpal)\bigr)\circ
exp\Bigl(Darit\bigl(\Delta(F)\bigr)\Bigr)\circ
Dgarit\bigl(\Delta^*(invpal)\bigr)^{-1},}$$
where the second equality follows from (3.3.4) (with $\frak{g}$, $exp_{\frak{g}}$
and $Ad_{\frak{g}}$ identified with $ARI_{ari}$, $exp_{ari}$ and
$Ad_{ari}$, and the same three terms for $\frak{h}$ with the corresponding
terms for $ARI_{Dari}$), the third from (3.2.4), the 
fourth from (3.3.2) and the fifth again from (3.2.4).
But the first and last expressions in (3.3.5) are equal to the exponentials 
of the left- and right-hand sides of (3.3.1). This concludes the proof of 
the Proposition.\hfill{$\diamondsuit$}
\vskip .5cm
We can now complete the proof of Proposition 2.2.1 by using Proposition 3.3.1 to 
compute the value
of $E_\psi(t_{02})$. By (3.2.9) and the Corollary to Proposition 3.2.3, we have
$$Dgarit\bigl(\Delta^*(invpal)\bigr)\cdot a=a-1+\Delta^*(invpal)
=ma\bigl(Ber_{-b}(a)\bigr)=ma(t_{02}).\eqno(3.3.6)$$
Recall that $E_\psi$ is nothing but the polynomial version of $Darit(M)$ 
restricted to the Lie algebra generated by the moulds $a$ and $B$.
Thus, to compute the value of $E_\psi$ on $t_{02}=Ber_{-b}(a)$, we can now
simply use (3.3.1) to compute the value of $Darit(M)$ on $ma(t_{02})$.
By (3.3.6), the rightmost map of the right-hand side of (3.3.1) maps $ma(t_{02})$ 
to $a$. By Lemma 3.1.1 (iv), the derivation $Darit(P)$ for any mould $P\in ARI$
extends to $a$ taking the value $P$ on $a$, so we can apply the middle 
map of (3.3.1) to $a$, obtaining
$$Darit\bigl(\Delta(F)\bigr)\cdot a=\Delta(F)
=dur\bigl(dar(F)\bigr)=ma\bigl([f(a,[b,a]),a]\bigr)$$
$$=ma\bigl([\psi(a,[a,b]),a]\bigr)=ma\bigl([\psi(a,t_{12}),a]\bigr).\eqno(3.3.7)$$
Finally, we note that by Lemma 3.2.2 (ii), the leftmost map of the right-hand side 
of (3.3.1) fixes $B_1=-ma(t_{12})$, so it also fixes $ma(t_{12})$.  By
(3.3.6), it sends $a$ to $ma(t_{02})$, so applying it to the rightmost term
of (3.3.7) we obtain the total expression
$$Darit(M)\bigl(ma(t_{02})\bigr)=ma\Bigl([\psi(t_{02},t_{12}),t_{02}]\Bigr).$$
In terms of polynomials, this gives the desired expression
$$E_\psi(t_{02})=[\psi(t_{02},t_{12}),t_{02}],$$
which concludes the proof.\hfill{$\diamondsuit$}

\vfill\eject
\noindent {\bf Appendix: Mould basics}
\vskip .3cm
For the purposes of this article, we use the term ``mould'' to refer only
to rational-function valued moulds with coefficients in $\Q$; thus, 
a mould is a family of functions $\{P_r(u_1,\ldots,u_r)\mid r\ge 0\}$ with 
$P_r(u_1,\ldots,u_r)\in \Q(u_1,\ldots,u_r)$.  In particular
$P_0(\emptyset)$ is a constant.  The {\it depth $r$} part of a mould is
the function $P_r(u_1,\ldots,u_r)$ in $r$ variables. By defining addition
and scalar multiplication of moulds in the obvious way, i.e.,~depth by
depth, moulds form a $\Q$-vector space that we call $Moulds$.  Following
Ecalle, we often drop the subscript $r$ from the mould notation;
i.e.,~we write $P(u_1,\ldots,u_r)$ to mean the rational function 
$P_r(u_1,\ldots,u_r)$, where the number of variables automatically indicates
which depth part we are considering.

We write $GARI$ for the set of moulds with $P(\emptyset)=1$, and
$ARI$ for the set of moulds\footnote{$^8$}{Ecalle
uses the notation $ARI$ for the space of these moulds equipped with the
$ari$-bracket, that we denote $ARI_{ari}$, and in fact he considers more
general {\it bimoulds} in two sets of variables.} with $P(\emptyset)=0$.
Then $ARI$ forms a vector subspace of $Moulds$. 

Let $(Moulds)^{pol}$ denote the subspace of polynomial-valued moulds,
i.e.,~moulds such that $P(u_1,\ldots,u_r)$ is a polynomial in each depth 
$r\ge 0$, and $ARI^{pol}$ the polynomial-valued subspace of $ARI$.
In this appendix we will stress the connections between polynomial-valued
moulds and power series in the non-commutative variables $a$ and $b$, 
showing in particular how familiar notions from multizeta theory (the 
Poisson-Ihara bracket, the twisted Magnus group etc.) not only translate 
over to the corresponding moulds, but generalize to all moulds.
\vskip .2cm
Let $C_i=ad(a)^{i-1}(b)$ for $i\ge 1$.  Let the depth of a monomial
$C_{i_1}\cdots C_{i_r}$ be the number $r$ of $C_i$ in the monomial;
the depth forms a grading on the free associative ring of polynomials
in the $C_i$'s. Let
$\Q\langle C\rangle =\Q\langle C_1,C_2,\ldots\rangle$
denote the depth completion of this ring, i.e., 
$\Q\langle C\rangle$ is the space of power series that are polynomials 
in each depth. We also write
$$L[C]={\rm Lie}[C_1,C_2,\ldots]\eqno(A.1)$$  
for the corresponding free Lie algebra. Note that the freeness follows 
from Lazard elimination, which also shows that we have the 
isomorphism
$$\Q a\oplus L[C]\simeq {\rm Lie}[a,b].$$
Ecalle uses the notation $ma$ to denote the standard vector space isomorphism from 
$\Q\langle C\rangle$ to the space 
$(Moulds)^{pol}$ of polynomial-valued moulds defined by
$$\eqalign{ma:\Q\langle C\rangle&\buildrel\sim\over\rightarrow (Moulds)^{pol}\cr
C_{k_1}\cdots C_{k_r}&\mapsto (-1)^{k_1+\cdots+k_r-r}
u_1^{k_1-1}\cdots u_r^{k_r-1}}\eqno(A.2) $$ 
on monomials and extended by linearity. This map $ma$ can also be 
considered as a ring isomorphism when $(Moulds)^{pol}$ is equipped with
the suitable multiplication, cf.~the remarks following (A.4) below.
(We use the same notation $ma$ when $C_i=ad(x)^{i-1}(y)$, for polynomials 
usually considered in ${\rm Lie}[x,y]$, such as polynomials in $\grt$.) 
For any map $\Phi:\Q\langle C\rangle\rightarrow \Q\langle C\rangle$, we
define its transport $ma(\Phi)$ to $(Moulds)^{pol}$, 
namely the corresponding map on polynomial-valued moulds
$$ma(\Phi):(Moulds)^{pol}\rightarrow (Moulds)^{pol}$$ 
by the obvious relation
$$ma(\Phi)(ma(f))=ma\bigl(\Phi(f)\bigr)\ \ {\rm for\ all\ }f\in \Q\langle C
\rangle.\eqno(A.3)$$

\noindent {\bf Power series, moulds, standard multiplication and Lie bracket.}
Via the map (A.2), many of the familiar notions 
associated with power series and Lie series pass to polynomial moulds, 
with general expressions that are in fact valid for all moulds.

In particular, the standard mould multiplication $mu$ is given by
$$mu(P,Q)(u_1,\ldots,u_r)=\sum_{i=0}^r P(u_1,\ldots,u_i)Q(u_{i+1},\ldots,u_r).$$
For simplicity, we write $P\,Q=mu(P,Q)$. The multiplication $mu$ 
generalizes ordinary
multiplication of non-commutative power series in the sense that
$$ma(fg)=mu\bigl(ma(f),ma(g)\bigr)=ma(f)\,ma(g)\eqno(A.4)$$
for $f,g\in \Q\langle C\rangle$.  The space $(Moulds)^{pol}$ is a ring
under the $mu$ multiplication, generated by the depth $1$ polynomial
moulds $B_i$ given by $B_i(u_1)=u_1^i$ for $i\ge 0$.  By (A.4), the linear 
map $ma$ from (A.2) can be defined as a ring isomorphism from
$\Q\langle C\rangle$ to $(Moulds)^{pol}$, taking values
$ma(C_i)=(-1)^{i-1}B_{i-1}$ on the generators $C_i$ for $i\ge 1$.

A mould $P$ is invertible for the $mu$-multiplication if and only if its
constant term $P(\emptyset)\in \Q$ is invertible.  If the constant term
is 1, the formula for the $mu$-inverse $P^{-1}=invmu(P)$ is explicitly
given by 
$$P^{-1}({\bf u})=\sum_{0\le s\le r}(-1)^s 
\sum_{{\bf u}={\bf u}_1\cdots {\bf u}_s}
P({\bf u}_1)\cdots P({\bf u}_s),$$
where the sum runs over all ways ${\bf u}_1\cdots {\bf u}_s$ of cutting
the word ${\bf u}=(u_1,\ldots,u_r)$ into $s$ non-empty chunks.
By (A.4), if $f\in\Q\langle C\rangle$ is invertible (i.e.,~has
non-zero constant term), we have $ma(f^{-1})=P^{-1}$.

The $mu$-multiplication makes $GARI$, the set of moulds with constant term 1,
into a group that we denote by
$GARI_{mu}$.  Defining the associated $lu$-bracket by $lu(P,Q)=mu(P,Q)-mu(Q,P)$,
i.e.,  $[P,Q]=P\,Q-Q\,P$, gives ARI the structure of a Lie algebra that we call 
$ARI_{lu}$.
\vskip .2cm
\noindent {\bf Mould symmetries.} A mould $P$ is said to be {\it alternal} if
$$\sum_{{\bf u}\in sh\bigl((u_1,\ldots,u_i),(u_{i+1},\ldots,u_r)\bigr)}
P({\bf u})=0\eqno(A.5)$$
for $1\le i\le r-1$. 

It is well-known that $p\in\Q\langle C\rangle$ satisfies the shuffle relations 
if and only if $p$ is a Lie polynomial, i.e., $p\in {\rm Lie}[C]$.  
The alternality property on moulds is analogous to these shuffle
relations, i.e.,~a polynomial $p\in\Q\langle C\rangle$ satisfies the shuffle 
relations if and only if $ma(p)$ is alternal.  
(See [S, \S 2.3 and Lemma 3.4.1.].)  This
shows that, writing $ARI^{al}$ for the subspace of alternal moulds and 
$ARI^{pol,al}$ for the subspace of alternal polynomial-valued moulds,  
the map $ma$ restricts to a Lie algebra isomorphism 
$$ma:{\rm Lie}[C]\buildrel{ma}\over\longrightarrow ARI_{lu}^{pol,al}.$$ 

Let the {\it swap} operator on moulds be defined by
$$swap(A)(v_1,\ldots,v_r)=A(v_r,v_{r-1}-v_r,\ldots,v_1-v_2).$$
Here the use of the alphabet $v_1,v_2,\ldots$ instead of $u_1,\ldots,u_r$
is purely a convenient way to distinguish a mould from its swap. The
mould $swap(A)$ is alternal if it satisfies the property (A.5) in the $v_i$'s.
The space of moulds that are alternal and have a swap that is also alternal
is denoted $ARI^{al/al}$; these moulds are said to be {\it strictly
bialternal}.  We particularly consider the situation where
a mould is alternal and its swap differs from an alternal mould by addition
of a constant-valued mould. Such moulds are called {\it bialternal}, and
the space of bialternal moulds is denoted $ARI^{al*al}$.  The space
of polynomial-valued bialternal moulds is denoted $ARI^{pol,al*al}$.
Finally, we recall that Ecalle uses the notation of underlining the symmetry
of a mould to indicate that its depth 1 part is an even function of $u_1$;
thus we use the notation $ARI^{pol,\underline{al}*\underline{al}}$ etc. to
denote the subspaces of moulds that are even in depth 1.  The subspace
$ARI_{ari}^{pol,\underline{al}*\underline{al}}$ forms a Lie algebra under the 
$ari$-bracket (cf. [S, Theorem 2.5.6]), which is isomorphic under the map 
$ma$ to the ``linearized double shuffle'' Lie algebra $\ls$ studied for 
example in [Br3].

Ecalle introduces a second symmetry called {\it alternility} on moulds
in the $v_i$'s, which generalizes the usual stuffle relations on polynomials
in $a$ and $b$.  As above, we write $ARI^{al/il}$, $ARI^{al*il}$ and
$ARI^{\underline{al}*\underline{il}}$ for the space of alternal moulds 
with swap that is alternil, resp. alternil up to addition of a constant
mould, resp. also even in depth 1.  The space
$ARI^{pol,\underline{al}*\underline{il}}$ is isomorphic under the map $ma$
to the double shuffle Lie algebra $\ds$. [S, ??]
\vskip .2cm
\noindent {\bf Twisted Magnus automorphism and group law.}
Let $G\subset \Q\langle C\rangle$ denote the set of
power series with constant term 1, so that $ma$ gives a bijection
$G\rightarrow GARI^{pol}$ to the set of polynomial-valued moulds with constant 
term 1.  We write ${\bf G}$ for the group obtained by putting
the standard power series multiplication on $G$, so that we have a group
isomorphism ${\bf G}\simeq GARI^{pol}_{mu}$. For all $p\in G$, 
we define the associated ``twisted Magnus'' automorphism $A_p$ of ${\bf G}$,
defined by $A_p(a)=a$, $A_p(b)=pbp^{-1}$. These automorphisms satisfy the 
composition law
$$(A_q\circ A_p)(b)=A_q(p)qbq^{-1}A(p)^{-1},$$
which defines a different multiplication on the set $G$, given by
$$p\odot q=A_q(p)q = p\bigl(a,qbq^{-1}\bigr)\,q(a,b),\eqno(A.6)$$
satisfying
$$A_{p\odot q}=A_q\circ A_p.$$
The inverse of the automorphism $A_p$ is given by $A_q$ where $q$ is
the unique power series such that the right-hand side of (A.6) is equal
to 1.  We write $G_\odot$ for the ``twisted Magnus'' group obtained by putting 
the multiplication law (A.6) on $G$.  The association $p\mapsto A_p$ extends to the 
general case of moulds by associating to every $P\in GARI$ the automorphism 
of $GARI_{mu}$ defined by Ecalle and denoted $garit(P)$, whose action
on $Q\in GARI$ is given by
$$\bigl(garit(P)\cdot Q\bigr)({\bf u})=\sum_{s\ge 0} 
\sum_{{\bf u}={\bf a}_1{\bf b}_1{\bf c}_1\cdots {\bf a}_s{\bf b}_s
{\bf c}_s} \!\!\!Q(\lceil {\bf b}_1\rceil\cdots \lceil {\bf b}_2\rceil)
P({\bf a}_1)\cdots P({\bf a}_s) P^{-1}({\bf c}_1)\cdots P^{-1}({\bf c}_s),$$
where the sum runs over all ways of cutting the word
${\bf u}=(u_1,\ldots,u_r)$ into $3s$ chunks of which the ${\bf b}_i$'s may
not be empty, ${\bf a}_1$ and ${\bf c}_s$ may be empty, and the
interior chunks ${\bf a}_i$ and ${\bf c}_j$ may be empty as long as no
interior double chunk ${\bf c}_i{\bf a}_{i+1}$ is empty.  Note that
because $GARI_{mu}$ is a huge group containing all possible moulds with
constant term 1, the automorphism $garit(P)$ cannot be determined simply 
by giving its value on some simple generators as we do for $A_p$.  
However, $garit(P)$ extends to $a$ taking the value $a$, and restricted
to the Lie algebra $(ARI^a_{lu})^{pol}$ generated by $a$ and $B$ (isomorphic 
to ${\rm Lie}[a,b]$), we find
$$garit(P)\cdot a=a,\ \ \ garit(P)\cdot B=PBP^{-1}.\eqno(A.7)$$

In analogy with 
the formula for $\odot$ given in (A.6), $garit$ defines a multiplication law $gari$ on $GARI$ by the formula
$$gari(P,Q)=mu\bigl(garit(Q)\cdot P,Q)=\bigl(garit(Q)\cdot P\bigr)\,Q.$$
We write $GARI_{gari}$ for the group obtained by equipping $GARI$ with this
multiplication.
\vskip .1cm
\noindent {\bf Poisson-Ihara bracket, exponential, linearization.}
For all $P\in ARI$, Ecalle defines a derivation $arit(P)$ of $ARI_{lu}$
by the formula
$$arit(F)\cdot M({\bf u})=\sum_{{\bf u}={\bf abc},{\bf c}\ne\emptyset} 
M({\bf a}\lceil{\bf c}) F({\bf b})-\sum_{{\bf u}={\bf abc},{\bf a}\ne \emptyset}
M({\bf a}\rceil{\bf c})F({\bf b}).$$
For $B=ma(b)$, i.e., $B(u_1)=1$, this formula yields
$$arit(P)\cdot B=[P,B].\eqno(A.8)$$

If $P=ma(f)$ for a polynomial $f\in {\rm Lie}[C]$, then $arit(P)$ restricts
to $ARI_{lu}^{pol,al}$, and as we saw in Lemma 3.1.1 (iii), it extends
to all of $(ARI^a_{lu})^{pol,al}$ taking the value $0$ on $a$.  It 
corresponds on the isomorphic Lie algebra ${\rm Lie}[a,b]$ to the Ihara 
derivation $D_f$ defined by
$$D_f(a)=0,\ \ D_f(b)=[f,b].\eqno(A.9)$$
The Lie bracket $\{\cdot,\cdot\}$
that we put on $L[C]$, known as the Poisson bracket or Ihara bracket,
comes from bracketing the derivations $D_f$, i.e.,
$$[D_f,D_g]=D_{\{f,g\}}\ \ {\rm where}\ \ \{f,g\}=D_f(g)-D_g(f)-[f,g].
\eqno(A.10)$$
We obtain a pre-Lie law by linearizing the multiplication law $\odot$
defined in (A.6).  In fact, because $\odot$ is linear in $p$, we only
need to linearize $q$, so we write $q=1+tf$ and compute the 
coefficient of $t$ in 
$$p\Bigl(a,(1+tf)b(1-tf)\Bigr)\bigl(1+tf(a,b)\bigr)=
p\bigl(a,b+t[f,b]\bigr)\bigl(1+tf(a,b)\bigr),$$
obtaining the expression 
$$p\odot f=pf+D_f(p),\eqno(A.11)$$
valid for all $p\in \Q\langle C\rangle$, $f\in L[C]$.
In particular, the pre-Lie law gives another, equivalent way to obtain the 
Poisson bracket, namely $\{p,q\}=p\odot q-q\odot p$.
The exponential map
$exp_\odot:L[C]_{\{\cdot,\cdot\}}\rightarrow G_\odot$ is
then defined via the pre-Lie law by
$$exp_\odot(f)=\sum_{n\ge 0}{{1}\over{n!}}f^{\odot n},\eqno(A.12)$$
where the pre-Lie law is composed from left to right, so that the rightmost
argument is always $f\in L[C]$.  The exponential map defined this way 
satisfies the basic identities
$$exp\bigl(D_f\bigr)=A_{exp_\odot(f)},\eqno(A.13)$$
and
$$exp(D_f)\circ exp(D_g)=exp\bigl(ch_{\{\cdot,\cdot\}}(D_f,D_g)\bigr),\eqno(A.14)$$
where $ch_{\{\cdot,\cdot\}}$ denotes the Campbell-Hausdorff law on $L[C]$ equipped with the Poisson-Ihara Lie bracket (A.10).
\vskip .2cm
All these standard constructions extend to the case of general moulds;
Ecalle gives explicit formulas for the pre-Lie law $preari$ and 
for the exponential $exp_{ari}$, namely
$$preari(P,Q)=PQ+arit(Q)\cdot P\ \ {\rm and}\ \
exp_{ari}(P)=\sum_{n\ge 0}{{1}\over{n!}} preari(\underbrace{P,\ldots,P}_n),$$
which clearly extend (A.11) and (A.12) above, and satisfy the analogous
formulas generalizing (A.13) and (A.14), namely
$$exp\bigl(arit(P)\bigr)=garit\bigl(exp_{ari}(P)\bigr)\eqno(A.15)$$
and
$$exp\bigl(arit(P)\bigr)\circ exp\bigl(arit(Q)\bigr)=exp\Bigl(ch\bigl(
arit(P),arit(Q)\bigr)\Bigr).\eqno(A.16)$$
The exponential maps satisfy the properties
$$exp\bigl(arit(P)\bigr)\circ exp\bigl(arit(Q)\bigr)=
exp\Bigl(arit\bigl(ch_{ari}(P,Q)\bigr)\Bigr) \eqno(A.17)$$
for the Campbell-Hausdorff law $ch_{ari}$ on $ARI_{ari}$. These
properties are expressed by the commutative diagram
$$\xymatrix{ARI\ar[r]^{exp_{ari}}\ar[d]_{arit}&GARI\ar[d]^{garit}\\
{\rm \Der}\,ARI_{lu}\ar[r]^{exp}&{\rm Aut}\,ARI_{lu}.}\eqno(A.18)$$
\vskip .1cm
We conclude this appendix with a linearization lemma used
in the proof of Proposition 3.2.3.
\vskip .3cm
\noindent {\bf Lemma A.1.} {\it Let $P\in ARI$. Then the derivation
$-arit(P)+ad(P)$ extends to $a$ taking the value
$[P,a]$ on $a$, and we have 
$$exp\bigl(-arit(P)+ad(P)\bigr)\cdot a= R^{-1}aR$$
where $R=exp_{ari}(-P)$.}
\vskip .2cm
\noindent Proof. Since $arit(P)$ extends to $a$ taking the value $0$ by
Lemma 3.1.1 (iii), it suffices to check that $ad(P)$ extends to $a$
via $ad(P)\cdot a = [P,a]$, i.e., that this action respects the formula 
$[Q,a]=dur(Q)$. Indeed, we have
$$ad(P)\cdot [Q,a]=[ad(P)\cdot Q,a]+[Q,ad(P)\cdot a]=
\bigl[[P,Q],a\bigr]+\bigl[Q,[P,a]\bigr]=
[P,[Q,a]]=ad(P)\cdot dur(Q).$$
For a real parameter $t\in [0,1]$, let $R_t=exp_{ari}(-tP)$,
and let $A_t$ denote the automorphism of $(ARI^a_{lu})^{pol}$ defined by
$$A_t(a)=R_t^{-1}aR_t, \ \ \ A_t(B)=B,$$
so that $A_1(a)=R^{-1}aR$. Let $D=log(A)$; we will prove that 
$D=-arit(P)+ad(P)$ on $(ARI^a_{lu})^{pol}$.  We compute $D(a)$
and $D(b)$ by the linearization formula
$$D(a)={{d}\over{dt}}|_{t=0}\bigl(A_t(a)\bigr)\ \ \ {\rm and}\ \ \
D(b)={{d}\over{dt}}|_{t=0}\bigl(A_t(b)\bigr).$$
The second equality yields $D(b)=0$.  Let us compute $D(a)$.  Using
$R_0=1$ and ${{d}\over{dt}}|_{t=0}R_t=-P$, we find
$$\eqalign{D(a)&={{d}\over{dt}}|_{t=0}\Bigl(A_t(a)\Bigr)\cr
&={{d}\over{dt}}|_{t=0}\Bigl(R_t^{-1}aR_t\Bigr)\cr
&=\Bigl(-R_t^{-1}{{d}\over{dt}}(R_t)R_t^{-1}aR_t+R_t^{-1}a{{d}\over{dt}} (R_t)
\Bigr)|_{t=0}\cr
&=Pa-aP.}$$
Thus $D(a)=[P,a]=\bigl(-arit(P)+ad(P)\bigr)\cdot a$ and $D(b)=0=
\bigl(-arit(P)+ad(P)\bigr)\cdot b$, which concludes the proof.\hfill{$\diamondsuit$}
\vskip 1cm
\noindent {\bf References}
\vskip .5cm
\noindent [B] S. Baumard, Aspects modulaires et elliptiques des relations
entre multiz\^etas, Doctoral thesis, 2014.
\vskip .1cm
\noindent [BS] S. Baumard, L. Schneps, On the derivation representation
of the fundamental Lie algebra of mixed elliptic motives, to appear
in {\it Annales Math\'ematiques du Qu\'ebec}, 2016.
\vskip .1cm
\noindent [Br1] F. Brown, Mixed Tate motives over ${\Bbb Z}$, {\it Ann. of
Math.} (2) {\bf 175}, no. 2 (2012), 949-976.
\vskip .1cm
\noindent [Br2] F. Brown, Talk on Anatomy of Associators,
online lecture notes, http://www.ihes.fr/~brown, 2013.
\vskip .1cm
\noindent [Br3] F. Brown, Zeta elements in depth 3 and the fundamental Lie 
algebra of the infinitesimal Tate curve, {\it Forum Math. Sigma} {\bf 5} (2017).
\vskip .1cm
\noindent [Br4] F. Brown, Depth-graded motivic multiple zeta values,
{\it Compos. Math.} {\bf 157} 3 (2021), 529-572.
\vskip .1cm
\noindent [E1] J. Ecalle, The flexion structure of dimorphy: flexion units, singulators,
generators, and the enumeration of multizeta irreducibles, in {\it Asymptotics
in Dynamics, Geometry and PDEs; Generalized Borel Summation II}, O. Costin,
F. Fauvet, F. Menous, D. Sauzin, eds., Edizioni della Normale, Pisa, 2011.
\vskip .1cm
\noindent [E2] J. Ecalle, Eupolars and their bialternality grid, 
{\it Acta Math. Vietnam.} {\bf 40} no. 4 (2015), 545-636. 
\vskip .1cm
\noindent [CEE] D. Calaque, B. Enriquez, P. Etingof, Universal KZB equations: the elliptic case, {\it Algebra, Arithmetic and Geometry: in honor of Yu.~I.~Manin}, Vol. I, 165-266, {\it Progr. Math.} {\bf 269}, Birkh\"auser Boston, 2009.
\vskip .1cm
\noindent [En1] B. Enriquez, Elliptic Associators, {\it Selecta Math.} {\bf 20} no. 2 
(2014), 491-584.
\vskip .1cm
\noindent [En2] B. Enriquez, Analogues elliptiques des nombres multizetas,
{\it Bull. Soc. Math. France} {\bf 144} no. 3 (2016), 395-427.
\vskip .1cm
\noindent [F] H. Furusho, Double shuffle relation for associators,
{\it Ann. of Math.} {\bf 174} no. 1 (2011), 341-360.
\vskip .1cm
\noindent [G] A. Goncharov, Galois symmetries of fundamental groupoids
and noncommutative geometry, {\it Duke Math J.} {\bf 128} (2005) no. 2,
209-284.
\vskip .1cm
\noindent [HM] R. Hain, M. Matsumoto, Universal mixed elliptic motives,
{\it J. Inst. Math. Jussieu} {\bf 19} no. 3 (2020), 663-766. 
\vskip .1cm
\noindent [LMS] P. Lochak, N. Matthes, L. Schneps, Elliptic multiple zeta 
values and the elliptic double shuffle relations, {\it IMRN} {\bf 2021} no. 1, 
695-753.
\vskip .1cm
\noindent [NTU] H. Nakamura, H. Tsunogai, R. Ueno, Some stability properties of
Teichm\"uller modular function fields with pro-$\ell$ weight structures,
{\it Math. Ann.} {\bf 302}, 197-213 (1995).
\vskip .1cm
\noindent [P] A. Pollack, Relations between derivations arising from modular
forms, Duke University Senior Thesis, 2009.
\vskip .1cm
\noindent [R] G. Racinet, S\'eries g\'en\'eratrices non commutatives de polyz\^etas
et associateurs de Drinfel'd, Doctoral thesis, 2000.
\vskip .1cm
\noindent [S] L. Schneps, ARI, GARI, Zig and Zag, An Introduction to Ecalle's 
theory of moulds, arXiv:1507.01534, 2015.
\vskip .1cm
\noindent [S2] L. Schneps, Double shuffle and Kashiwara-Vergne Lie algebras,
{\it J. Algebra} {\bf 367} (2012), 54-74.
\vskip .1cm
\noindent [So] I. Soud\`eres, Motivic double shuffle, {\it Int. J.
Number Theory} {\bf 6} (2010) no. 2, 339-370.
\vskip .1cm
\noindent [T1] H. Tsunogai, On some derivations of Lie algebras related
to Galois representations, {\it Publ. RIMS} {\bf 31} No. 1 (1995), 113-134.
\vskip .1cm
\noindent [T2] H. Tsunogai, The stable derivation algebras for higher
genera, {\it Israel J. Math.} {\bf 136} (2003), 221-250.
\bye